\documentclass{amsart}
\usepackage{amsmath,amssymb}
\usepackage{amscd}
\usepackage[english]{babel}

\usepackage{stmaryrd}

\newcommand{\R}{\mathbb{R}}
\newcommand{\N}{\mathbb{N}}

\newcommand{\with}{\;:\;}
\newcommand{\id}{\hbox{id}}

\usepackage[utf8]{inputenc}  
\usepackage[T1]{fontenc}

\newtheorem{theorem}{Theorem}[section]
\newtheorem{proposition}[theorem]{Proposition}
\newtheorem{lemma}[theorem]{Lemma}
\newtheorem{corollary}[theorem]{Corollary}

\theoremstyle{definition}
\newtheorem{definition}[theorem]{Definition}
\newtheorem{remark}[theorem]{Remark}

\begin{document}

\markboth{M. Bessa \& M. Stadlbauer}{On the Lyapunov spectrum of relative transfer operators}


\title{On the Lyapunov spectrum of relative transfer operators}

\author{M\'{A}RIO BESSA}

\address{
Departamento de Matem\'atica,  Universidade da Beira Interior\\ Rua Marqu\^es d'\'Avila e Bolama,
6201-001 Covilh\~a, Portugal\\
bessa@ubi.pt
}

\author{MANUEL STADLBAUER}

\address{
Departamento de Matem\'atica, Universidade Federal do Rio de Janeiro\\
C. P. 68.530, 21941-909 Rio de Janeiro (RJ), Brazil.\\
manuel@im.ufrj.br
}

\maketitle


\begin{abstract}
We analyze the Lyapunov spectrum of the relative Ruelle operator associated with a skew product whose base is an ergodic automorphism and whose fibers are full shifts. We prove that these operators can be approximated in the  $C^0$-topology by positive matrices with an associated dominated splitting. 
\end{abstract}

\section{Introduction: Basic definitions and statement of the results}

We begin with the definition of the class of the relevant dynamical systems and cocycles. For a fixed $n \in \N$, let $\Sigma :=  \{1, 2,  \ldots, n\}^\N$  refer to the one-sided full shift space with $n$ states and $\theta: \Sigma \to \Sigma, (w_0 w_1\ldots )\mapsto (w_1\ldots )$ be the left shift. Furthermore, for a compact Hausdorff space $M$ and a  continuous, ergodic automorphism  $f\colon M\rightarrow M$   
with respect to the $f$-invariant Borel probability  $\mu$ on $M$, let 
\[ T:\Sigma \times M \to \Sigma \times M, (w ,x)\mapsto (\theta(w),f(x)).  \]
We are  interested in the time evolution of the family  $\{\mathcal{L}_x: x \in M \}$ of Ruelle operators  acting on H\"older functions defined by  
\[  \mathcal{L}_x(h)(w):=  \sum_{v: \theta(v)=w}  h(v) e^{\varphi_x(v)}, \]
where $\varphi : \Sigma \times M \to \R$, $(w,x) \mapsto \varphi_x(w)$ is a \emph{uniformly relatively H\"older continuous}
 potential function, that is  $x \mapsto  \varphi_x(w)$ is continuous and $w \mapsto \varphi_x(w)$ is H\"older continuous with continuously varying H\"older coefficients (see Def. \ref{def:relatively_Hoelder}).
As usual, the time evolution is defined by
\[ 
\mathcal{L}_x^k := \left\{ \begin{array}{ll}
\mathcal{L}_{f^{k-1}(x)}\cdots \mathcal{L}_x & \hbox{if } k>0 \\
\id & \hbox{if } k=0, \end{array}\right. \]
and, in particular, the {cocycle identity} $\mathcal{L}^{k+n}_x = \mathcal{L}^{k}_{f^n(x)} \circ \mathcal{L}^{n}_x$ holds. As it is well known from the theory of random dynamical systems, $\mathcal{L}_x$ acts on the space of H\"older functions and moreover, a fibrewise Perron-Frobenius-Ruelle theorem holds (see \cite{BG}). That is, there exist uniquely determined families $\{h_x : x \in M\}$ of strictly positive H\"older functions, $\{m_x : x \in M\}$ of probability measures and $\{\rho_x : x \in M\}$ of strictly positive constants such that $\mu$-a.s., with $\mathcal{L}_x^\ast$ referring to the dual of $\mathcal{L}_x$,
\[ \mathcal{L}_x(h) = \rho_x  h_{f(x)} ,\quad \mathcal{L}_x^\ast(m_{f(x)}) = \rho_x m_x. \]
If $\varphi$ is H\"older continuous as a function on $\Sigma \times M$, then $(x,w) \mapsto h_x(w)$ is globally H\"older continuous by application of the main result by Denker and Gordin in \cite{DenkerGordin:1999} in this restricted setting. 

The main intention of this article is to refine these results by showing that $\{\mathcal{L}_x\}$ can be approximated with respect to the uniform operator norm $\| \cdot \|_\infty$ by finite dimensional operators who admit a dominated splitting. In order to do so, we partially extend the result by Denker and Gordin to potentials who are only continuous in $x$ by proving uniform exponential decay for ratios of the iterated operators (see Theorem \ref{theo:ratio-exponential-decay} below). In here, $D_r(\cdot)$ refers to the $r$-H\"older coefficient and $\| \cdot \|_\mathcal{H}$ to the H\"older norm.

\medskip
\noindent \textbf{Theorem A.}
\textit{There exist $C > 0$ and $s \in (0,1)$ such that, for $g,h$ in the space of H\"older functions $\mathcal{H}$ and $g>0$, we have  that}
\[ \left\| \frac{\mathcal{L}^{n}_x(g) }{\mathcal{L}^{n}_x(h)}  - \frac{\int g dm_x}{\int h dm_x} \right\|_\mathcal{H}  \leq  C s^n \left(D_r(g) + \left| {\textstyle \frac{\int g dm_x}{\int h dm_x} }\right| D_r(h) \right) \|1/h\|_\infty. \]
Note that this is a result of independent interest since it does not require knowledge about $\{h_x\}$ and $\{\rho_x\}$. It is also worth noting that several aspects of the proof seem to be new. By considering a family of ratios of operators (as defined in (\ref{eq:ratio-operators})), it is possible to obtain exponential decay without a priori requiring the existence of $\{m_x\}$ or $\{h_x\}$ (see Theorem  \ref{label:main_prop_for_contraction}). As a corollary, we obtain existence, uniqueness  of the families $\{h_x\}$,  $\{m_x\}$ and $\{\rho_x\}$ (see Propositions \ref{prop:continuity_1} and \ref{prop:continuity_2}). 

We now recall the definitions of Lyapunov exponents, the Oseledets splitting and dominated splittings for linear cocycles. That is, for $m \in \N$ with  $\mathcal{M}_+^m$ referring to the set of positive $m\times m$ matrices, we consider the space of cocycles $C^0(M,\mathcal{M}_+^m)$ endowed with the topology defined by the norm
\[ \llbracket A\rrbracket_p :=\max_{x\in M} \sup_{v \in \R^m, \|v\|_p=1}\|A_x(v)\|_p,  \]
 where $\|\cdot\|_p$  refers to the $\ell_p$-norm on $\R^m$, for $p \in [1,\infty]$. If $p$ is not specified, then $p$ is considered to be equal to $2$. Note that the topology does not depend on the choice of $p$ as all norms on $\R^m$ are equivalent. As above, set 
$A^k_x:= A_{f^{k-1}(x)}\cdots A_x$ for  $ k>0$ and $A^0 = \id $.
Since $f$ leaves the probability $\mu$ invariant and $\mu$ is ergodic, the Oseledets theorem for non-invertible cocycles (see \cite{FroylandLloydQuas:2010}) guarantees that there exist $k\in \{1,\ldots,m\}$,  $\lambda_{1}>...>\lambda_{k} \geq -\infty$ (called \emph{Lyapunov exponents}), and  for $\mu$-a.e. point $x\in{M}$ a splitting $\mathbb{R}^m=E^{\lambda_1}_{x}\oplus...\oplus{E^{\lambda_k}_{x}}$ (called \emph{Oseledets splitting}) 
such that $A_x(E^{\lambda_i}_{x})=E^{\lambda_i}_{f(x)}$ if $\lambda_i >-\infty$, $A_x(E^{-\infty}_{x}) \subset E^{-\infty}_{f(x)}$, 
and, for  $v\in{E^{\lambda_i}_{x}\setminus\{0\}}$ and $i\in\{1,\ldots,k\}$,
\[
{\lim_{n \to \infty}} \; \frac{1}{n}\log{\|A^{n}_x v^{i}\|=\lambda_{i}}.
\]
Note that the Oseledets theorem does not make any statement on the continuity of the splitting. On the other hand, the notion of {projective hyperbolicity} is a statement about continuity and uniform separation of the Oseledets subspaces. That is, we say that a nontrivial $A$-invariant splitting $\mathbb{R}^m_{\Lambda}=F^1_{\Lambda}\oplus F^2_{\Lambda}$ over an $f$-invariant set $\Lambda\subset M$  is a $k$-\emph{dominated splitting} if $E^{-\infty}_\Lambda$ is well defined (possibly empty) and  a  subset of  $F^2_{\Lambda}$ and, for every $x\in \Lambda$,
\begin{equation}\label{ds}
 \frac{\|A^{k}_x|_{F^2_{x}}\|}{\mathfrak{m}(A^{k}_x|_{F^1_{x}})}\leq \frac{1}{2}.
\end{equation}
In here, $\mathfrak{m}$ refers to the co-norm of the operator, that is $\mathfrak{m}(A)=\|A^{-1}\|^{-1}$. Note that the co-norm in this situation is well defined since $E^{-\infty}_\Lambda \subset F^2_{\Lambda}$ implies that $\ker(A^{m}_x|_{F^1_{x}})=\{0\}$, and hence $A^{k}_x|_{F^1_{x}}$ is invertible. More generally, we call a $A$-invariant splitting $F^1 \oplus F^2 \oplus \cdots \oplus F^k$ \emph{dominated} if, for each $j = 1,\ldots,k-1$, the splitting $(F^1 \cdots F^j)\oplus(F^{j+1} \cdots F^k)$ is dominated.

The main result of the paper is a result on the approximation of $\{\mathcal{L}_x\}$ with respect to the uniform norm by finite dimensional operators admitting a dominated splitting.  In order to formulate the result, we refer to 
$E^{<\lambda}$ as the bundle of elements with Lyapunov exponent smaller than $\lambda$, that is, for $x \in M$, 
\[ E^{<\lambda}_x:= \{ v \in \R^n \with  \limsup_{n \to \infty} \frac{1}{n} \log \|A^{n}_x v\| < \lambda \}.\]

\medskip
\noindent \textbf{Theorem B.} \textit{Suppose that $\varphi$ is uniformly relatively H\"older continuous. 
Then, for each $\varepsilon >0$, there exist $m \in \N$, a continuous family of positive operators $\{A_x : \mathcal{H} \to \mathcal{H} \}$,  a continuous bundle $\{W_x : x \in  M\}$ of $m$-dimensional subspaces of $\mathcal{H}$ such that, for all $x \in M$,  
\begin{enumerate}
\item  $h_x \in W_x$ and $W_x$ is $\varepsilon$-dense in $\{ g \in \mathcal{H}:\|g\|_\mathcal{H}\leq 1\}$  with respect to $\|\cdot\|_\infty$, 
\item $A_x$ maps $W_x$ to $W_{f(x)}$ and $\|A_x(g) - \mathcal{L}_x(g)\|_\infty \leq \varepsilon \|g\|_\infty$ for all $g\in W_x$.
\end{enumerate}
Furthermore, the maximal Lyapunov exponent $\lambda_1$ associated with  $A$ is equal to $\int \log \rho_x d\mu$, $E_x^{\lambda_1}= \{ th_x :t \in \R\}$, the splitting $E^{\lambda_1}_M \oplus E^{<\lambda_1}_M \subset W_M$ is dominated, and either 
\begin{enumerate}
\item the Lyapunov spectrum of $A$ contains two points, or 
\item the Lyapunov spectrum of $A$ contains at least three points and the Oseledets subbundles define a dominated splitting.
\end{enumerate}}

Observe that the importance of the theorem stems from its relation to continuity of the Lyapunov exponents. In order to do so, recall that a dominated splitting automatically is continuous (see, e.g., § 2 in \cite{BV}). In particular, $E^{<\lambda_1}_M$ always is a continuous bundle and, if the Lyapunov spectrum contains at least three points, then also the Oseledets splitting is continuous. In particular, the assignment of the Lyapunov exponent to a point is globally defined and varies continuously with respect to the base point. However, if the spectrum only contains two points, then the second Lyapunov exponent might only exist almost everywhere. 

The main ingredients of the proof are Theorem A above and Theorem C below. Namely, $E_x^{\lambda_1}= \{ th_x :t \in \R\}$ and the dominated splitting property of $E^{\lambda_1}_M \oplus E^{<\lambda_1}_M$ are consequences of Theorem A. Furthermore, the continuity of the families $\{h_x\}$ and $\{\rho_x\}$ give rise to a continuous conjugation of $\{\mathcal{L}_x\}$ to a family of positive operators who leave invariant the constant function $\mathbf{1}$. Hence, in order to obtain the dichotomy concerning the continuity of Lyapunov exponents, it is necessary to make use of  the result of Bochi and Viana in \cite{BV} that a linear cocycle with values in an accessible group generically either has a trivial Lyapunov spectrum or the Oseledets splitting is dominated  (for accessibility and the Bochi-Viana result, see Definition \ref{def:accessible} and Lemma \ref{corollaryBV} below). We obtain the following result for the semigroup of stochastic matrices 
\[\mathcal{S}^n := \left\{ (a_{ij})\in \R^{n \times n} :  a_{ij} \in [0,1], \; 
 {\textstyle \sum_{j=1}^n a_{ij} =1 }\; \forall i = 1,2,\ldots,n \right\}.\]

\medskip
\noindent \textbf{Theorem C.}
\textit{  If $n \geq 2$, then there exists a residual subset $\mathcal{R}\subset C^0(M,\mathcal{S}^n)$ with the following properties. For all $S\in\mathcal{R}$, $\dim E^0_M =1$
 and the splitting $E^0_M \oplus E^{<0}_M$ is dominated. Furthermore, we either have that 
\begin{enumerate}
\item the Lyapunov spectrum of $S$ contains two points, or 
\item the Lyapunov spectrum of $S$ contains at least three points and the Oseledets subbundles define a dominated splitting.
\end{enumerate}}

The ideia of proof is to consider the induced action of a given $S$ on $E^{<0}_M$ and show that this action is close to an action with values in an accessible group. It is worth noting, that the main argument for this, Lemma \ref{lemma:projective_equality} below, seems to be optimal since the method of proof does not leave much flexibility for improvements. 

The article is structured as follows. In Section \ref{sec:rpf}, we give the details with respect to the topology of $\mathcal{H}$, introduce the Vaserstein metric and prove Theorem A by adaption of a result for non-stationary shift spaces in \cite{Stadlbauer:2015}. In Section \ref{sec:normal}, we consider stochastic cocycles and prove that 
 dominated splittings for the induced action can be lifted. Note that the dichotomy of Bochi and Viana cannot be applied immediately since $\mathcal{S}^n$ is not accessible (see Remark \ref{remark:non-accessible}).

We remark that Theorem B sheds light to the r\^{o}le of the relative Perron-Frobenius-Ruelle theorem and the dichotomy by Bochi and Viana in the context of transfer operators of non-invertible skew products. Namely, the relative Perron-Frobenius-Ruelle theorem guarantees that $E^{\lambda_1}_{x}$ is always one-dimensional whereas the result by Bochi and Viana provides an  approximation of the action of $\{\mathcal{L}_x\}$ on $E^{<\lambda_1}_M$ by dominated splittings.

\section{Exponential decay and continuity}\label{sec:rpf}
Throughout this section, if not stated explicitly, we do not to assume that $M$ is a topological space or that $f$ invertible. Continuity and invertibility only will be essential in the last part of this section in order to obtain continuity and uniqueness, respectively. 

We begin with the definition of the relevant maps, function spaces and operators. We refer to $\mathcal{W}^m := \{(a_1 \cdots a_m): 1 \leq a_i \leq m\}$ as the set of words of length $m$ and, for  $a=(a_1 \cdots a_m) \in \mathcal{W}^m$, the set 
\[ [a] = [a_1,...,a_m] :=\{  (w_i) \in \Sigma:\, w_i=a_i,\, i=1,...,m\}\]
is referred to as a cylinder and the map 
\[ \tau_a : \Sigma \to [a], (w_1 w_2 \cdots) \mapsto (a_1\cdots a_m w_1 \cdots) \]
 as the inverse branch defined by $a$. As it is well known, 
$d((v_i),(w_i)) := 2^{-\min \{i : v_i \neq w_i \}}$ defines a metric on $\Sigma$. Also recall that, with respect to this topology, $\Sigma$ is compact, $\theta$ is continuous, cylinder sets are clopen (closed and open) and $\tau_a $ is a homeomorphism. 
 Furthermore, for $g:\Sigma \to \R$ and $r \in \R$, we refer to  
 \[ D_r(g) := \sup \left\{ (g(v)-g(w))/r^{m} : v,w \in [a], a \in \mathcal{W}^m, m=1,2, \ldots    \right\} \]
 as the $r$-H\"older coefficient of $g$. The space of $r$-H\"older functions is then defined by 
 \[ \mathcal{H}_r := \{g : \|g \|_{\mathcal{H}} < \infty\}, \hbox{ with } \|g \|_{\mathcal{H}} := \|g\|_\infty + D_r(g).\]
We are now in position to specify the continuity assumptions on $\varphi_x$. That is, we say that $\varphi: \Sigma \times M \to \R, (w,x) \mapsto \varphi_x(w)$ has \emph{uniformly bounded H\"older coefficients} if there exists $r \in (0,1)$ such that $\sup_x D_r(\varphi_x) < \infty$.

Observe that this notion of H\"older continuity gives rise to the following bounded distortion estimate. That is, with $C_\varphi := \exp (\sup_x D_r(\varphi_x) /(1-r)$ and for $a \in \mathcal{W}^m$ and $v,w \in [a]$, it follows that 
\begin{equation} \label{eq:bounded_distortion}
C_\varphi^{-1} \leq  e^{\sum_{k=0}^{n-1} \varphi_{f^k(x)}(\theta^k(v)) - \sum_{k=0}^{n-1} \varphi_{f^k(x)}(\theta^k(w))} \leq  C_\varphi.  
 \end{equation} 
The proof of the estimate is well-known and therefore omitted (see, e.g. equation (3) in \cite{Stadlbauer:2015}). However, note that the estimate has various important consequences and implies, e.g., that $\mathcal{L}_{x}^{m}(\mathbf{1})(v)/\mathcal{L}_{x}^{m}(\mathbf{1})(w)\leq C_\varphi$ for all $v,w \in \Sigma$.  
 
Observe that, by a change to the equivalent metric $d((v_i),(w_i)) := r^{\min \{i : v_i \neq w_i \}}$, functions in $ \mathcal{H}_r$ are Lipschitz continuous. In particular, $\mathcal{H}_r$ is closely related to weak convergence
of measures and the Vaserstein distance through Kantorovich's duality.
Recall that, for two Borel probability measures $\nu_1,\nu_2$ on $\Sigma$, the set of couplings $\Pi(\nu_1,\nu_2)$ is defined as the set of Borel probability measures on $\Sigma \times \Sigma$ whose marginals are $\nu_1$ and $\nu_2$, respectively. The \emph{Vaserstein distance}  is  defined by 
\[W_r(\nu_1,\nu_2) := \inf\left\{ {\textstyle \int d_r(x,y) dQ} : {Q \in \Pi(\nu_1,\nu_2)} \right\}.\]
As it is well-known (see, e.g. \cite{Villani:2009}), this distance is compatible with weak convergence and Kantorovich's duality gives that 
\begin{equation}
\label{eq:KM-duality}
W_r(\nu_1,\nu_2) = \sup\left\{ {\textstyle \int g d\nu_1 - \int g d\nu_2} : D_r(g) \leq 1 \right\}.
\end{equation}
We now introduce the relevant operators in order to obtain a contraction of the Vaserstein distance. However, since the constant function $\mathbf{1}$ is not necessarily $\mathcal{L}_{x}$-invariant and the eigenfunctions $\{h_x\}$ are not yet known, we consider the operators, 
for $m,n \in \N$ and $g \in  \mathcal{H}_r$, 
\begin{equation}
\label{eq:ratio-operators}
\mathbb{P}^x_{m,n}(g):= \frac{\mathcal{L}_{f^n(x)}^m (g \cdot \mathcal{L}_{x}^{n}(\mathbf{1}))}{\mathcal{L}_{x}^{m+n}(\mathbf{1})} = \frac{\mathcal{L}_{x}^{m+n} (g\circ \theta^n)}{\mathcal{L}_{x}^{m+n}(\mathbf{1})}. 
\end{equation}
We now proceed by showing that this family of operators share many features known from Ruelle operators with normalised potentials. 
\begin{lemma}\label{lemma_doeblin_fortet_inequality}
 The operator $\mathbb{P}^x_{m,n}$ acts on $\mathcal{H}_r$. Furthermore, $\mathbb{P}^x_{m,n}(\mathbf{1}) = \mathbf{1}$ and $ D_r( \mathbb{P}^x_{m,n}(g)) \leq  C_\varphi  (2 \|g\|_\infty + D_r(g) r^m)$, for all $g \in \mathcal{H}_r$.
\end{lemma}
\begin{proof} For $v,w \in \Sigma$, we have  
\begin{eqnarray*}
\left|\mathbb{P}^x_{m,n}(g)(v)-\mathbb{P}^x_{m,n}(g)(w)\right| 
&\leq & \left| \frac{\mathcal{L}_{x}^{m+n} (g\circ \theta^n)(v) - \mathcal{L}_{x}^{m+n} (g\circ \theta^n)(w)}{\mathcal{L}_{x}^{m+n}(\mathbf{1})(v)}\right| \\ &+& \mathbb{P}^x_{m,n}(|g|)(w) 
\left| \frac{\mathcal{L}_{x}^{m+n} (\mathbf{1})(w) - \mathcal{L}_{x}^{m+n} (\mathbf{1})(v)}{\mathcal{L}_{x}^{m+n}(\mathbf{1})(v)}\right|.
\end{eqnarray*}
By applying a standard estimate (see, e.g., Prop. 2.1 in \cite{Stadlbauer:2015}), we have that
\begin{eqnarray*}
&& \left| {\mathcal{L}_{x}^{m+n} (g\circ \theta^n)(v)-\mathcal{L}_{x}^{m+n} (g\circ \theta^n)(w)}\right|\\ 
&\leq&
 \left( C_\varphi  \mathcal{L}_{x}^{m+n} (\mathbf{1})(v) \|g\circ \theta^n \|_\infty +    \mathcal{L}_{x}^{m+n} (\mathbf{1})(w) r^{m+n} D_r(g\circ \theta^n)\right)  {d_r(v,w)}\\
 &\leq&
\left( C_\varphi  \mathcal{L}_{x}^{m+n} (\mathbf{1})(v) \|g\|_\infty +    \mathcal{L}_{x}^{m+n} (\mathbf{1})(w) r^{m} D_r( g) \right)  {d_r(v,w)}.
\end{eqnarray*} 
Combining the above inequalities  with $\mathcal{L}_{x}^{m+n} (\mathbf{1})(v)/\mathcal{L}_{x}^{m+n} (\mathbf{1})(w) \leq C_\varphi$ and $\mathbb{P}^x_{m,n}(|g|)(w)  \leq \|g\|_\infty$
implies that $D_r( \mathbb{P}^x_{m,n}(g)) \leq  C_\varphi  (2 \|g\|_\infty + D_r(g) r^m)$. The remaining assertions are obvious.
\end{proof}

As consequence, we obtain from the above Lemma that the dual of $\mathbb{P}^x_{m,n}$ acts on the space of probability measure by $\int f d(\mathbb{P}^x_{m,n})^\ast(\nu) := \int \mathbb{P}^x_{m,n}(f) d\nu$. Following ideias in \cite{Hairer}, we now apply a further change of metric in order to obtain a contraction. In order to do so, set 
\[ \tilde{d}(v,w) := \min\{ 1, \alpha  d_r(v,w) \}, \hbox{ with  } \alpha:= 4C_\varphi\]
and let $\tilde{W}$ and $\tilde{D}$ refer to the corresponding Vaserstein distance and H\"older coefficient, respectively. Observe that  $d_r \leq \tilde{d} \leq \alpha d_r$ which implies that $W_r \leq \tilde{W} \leq \alpha W_r$ and $\tilde{D} \leq  D_r  \leq  \alpha \tilde{D}$. The following Theorem is an adaption of Lemma 2.1 in \cite{Stadlbauer:2015} (see also \cite{Kloeckner:2015aa}).

\begin{theorem} \label{label:main_prop_for_contraction} Suppose that $\varphi$ has uniformly bounded H\"older coefficients. Then there exist $k_0 \in \N$ and $s \in (0,1)$ such that, for all $n,m \in \N$ with $m \geq k_0$ and all Borel probability measures $\nu_1 , \nu_2$ and $g \in  \mathcal{H}_r $,
\begin{eqnarray}
\label{eq:contraction_vaserstein}
\tilde{W}((\mathbb{P}^x_{m,n})^\ast(\nu_1),(\mathbb{P}^x_{m,n})^\ast\ (\nu_2)) &\leq & s^{m} 
\tilde{W}( \nu_1 , \nu_2)\\
\label{eq:contraction_Lipschitz}
\tilde{D}(\mathbb{P}^x_{m,n})(g)  &\leq& s^{ m } \tilde{D}(g)
\end{eqnarray}
\end{theorem}

\begin{proof}  In the first three parts of the proof, we adapt the arguments in \cite{Stadlbauer:2015} to $\{\mathbb{P}^x_{m,n}\}$ and therefore omit the technical details in part (3), who allow to extend the result for Dirac measures to arbitrary probability measures.

\paragraph{(1) Local contraction.}
Assume that $\tilde{d}(v,w) < 1$, $v \neq w$ and $g \in \mathcal{H}_r$. Since $\mathbb{P}^x_{m,n}(\mathbf{1})=  \mathbf{1}$, we may assume without loss of generality for the estimate of $\tilde{D}(\mathbb{P}^x_{m,n}(g))$, that $ \inf_v g(v)=0$. Hence, for $m \geq -\log (\alpha)/\log(r)$, it follows from Lemma \ref{lemma_doeblin_fortet_inequality} that 
\begin{eqnarray*}
\frac{\left|\mathbb{P}^x_{m,n}(g)(v)-\mathbb{P}^x_{m,n}(g)(w)\right| }{ \tilde{d}(v,w) }
&\leq &
	 \frac{2C_\varphi \|g\|_\infty + C_\varphi D_r(g) r^m }{\alpha} \\
& \leq & 
	\frac{\|g\|_\infty}{2} + { C_\varphi \tilde{D}(g) r^m} \leq \frac{\tilde{D}(g)}{2} +  \frac{\tilde{D}(g)}{4}. 
\end{eqnarray*}
With $\delta_v$ referring to the Dirac measure in $v$, we have $\int g d(\mathbb{P}^x_{m,n})^\ast(\delta_v) = \mathbb{P}^x_{m,n}(g)(v)$ and $\tilde{W}( \delta_v, \delta_w)=\tilde{d}(v,w)$. Hence, by (\ref{eq:KM-duality}), it follows that 
\[ \tilde{W}((\mathbb{P}^x_{m,n})^\ast(\delta_v),(\mathbb{P}^x_{m,n})^\ast\ (\delta_w)) \leq  \frac{3}{4} \tilde{W}( \delta_v, \delta_w).\]
\paragraph{(2) Global contraction.} Assume that $\tilde{d}(v,w) = 1$ with  $v \neq w$. Furthermore, for $m \in \N$ and $a \in \mathcal{W}^m$, set 
$\Phi_a(v):= \exp({\sum_{k=0}^{m-1}  \varphi_{f^k(x)}(\theta^k(\tau_a(v)))}) $
and define a measure $R^m_{v,w}$ on $\Sigma \times \Sigma$ by 
\[R^m_{v,w} = \frac{\sum_{a \in \mathcal{W}^m} \min \{ \Phi_a(v),\Phi_a(w) \} \delta_{(\tau_a(v),\tau_a(w)}}{\max\{ \mathcal{L}^m_x(\mathbf{1})(v),\mathcal{L}^m_x(\mathbf{1})(w)\}}. \]
Note that  $R^m_{v,w} (\Sigma \times \Sigma)\leq 1$. Moreover, by combining (\ref{eq:bounded_distortion}) with $\tilde{d}(\tau_a(v),\tau_a(w)) \leq \alpha r^m$, we have for $\Delta_m := \{ (u_1,u_2): \tilde{d}(u_1,u_2)  \leq \alpha r^m \}$ that $R^m_{v,w} (\Delta_m) \geq \frac{1}{C_\varphi^2}$.
As it is possible to extend $R^m_{v,w}$ to an element in $\Pi(\delta_v,\delta_w)$, i.e. there exists a finite measure $Q$ such that  $Q^m_{v,w} = R^m_{v,w} + Q \in \Pi(\delta_v,\delta_w)$,
we obtain
\begin{eqnarray*}
\tilde{W}((\mathbb{P}^x_{m,n})^\ast(\delta_v),(\mathbb{P}^x_{m,n})^\ast\ (\delta_w)) &\leq& 
\int \tilde{d}(u_1,u_2) dQ^m_{v,w} \\
&\leq&  \alpha r^m Q^m_{v,w}(\Delta_m)  +  Q^m_{v,w}(\Delta_m^c) \\
&= & 1 - (1- \alpha r^m)  Q^m_{v,w}(\Delta_m) \leq 1 - \frac{1- \alpha r^m}{C_\varphi^2} .  
\end{eqnarray*}
Since $\tilde{d}(v,w) = 1$, we obtain a contraction, provided that $\alpha r^m <1$. 

\paragraph{(3) Combining (1) and (2).} For 
$k_0 >  -\log (\alpha)/\log(r)$ fixed and $s:=  \max(\{ 3/4, 1 - ({1- \alpha r^{k_0}})/{C_\varphi^2}\}$, parts (1) and (2) imply that 
\begin{equation}\label{eq:contraction_Dirac} \tilde{W}((\mathbb{P}^x_{k,n})^\ast(\delta_v),(\mathbb{P}^x_{k,n})^\ast\ (\delta_w)) \leq  s \tilde{W}( \delta_v, \delta_w), \hbox{ for all } k \geq k_0.\end{equation}
Hence, by Kantorovich's duality, $\tilde{D} \circ \mathbb{P}^x_{k,n} \leq s \tilde{D}$. It then follows either from general ideias from optimal transport or from the density of Dirac mesures that (\ref{eq:contraction_Dirac})
holds for all Borel probability measures. The proofs can be found in \cite{Kloeckner:2015aa} or \cite{Stadlbauer:2015} and are therefore omitted. Hence, we have shown that, if $k \geq k_0$, then
\[ \tilde{W}((\mathbb{P}^x_{k,n})^\ast(\,\cdot\,), (\mathbb{P}^x_{k,n})^\ast(\, \cdot\,)) 
\leq s 
\tilde{W}(\,\cdot\,,\,\cdot\,).\]
\paragraph{(4) Iterates.} First observe that 
$\mathbb{P}^x_{m,n} = \mathbb{P}^x_{m-j,n+j} \circ \mathbb{P}^x_{j,n}$, for $j,m,n \in \N$ with $m>j$. 
By induction, this implies that
\begin{eqnarray*}
 \mathbb{P}^x_{kl + j,n}
 & = &  \mathbb{P}^x_{kl,n + j} \circ   \mathbb{P}^x_{j,n} =  
 \mathbb{P}^x_{k(l-1),n+k + j} \circ \mathbb{P}^x_{k,n + j} \circ   \mathbb{P}^x_{j,n} \\
 & = & \mathbb{P}^x_{k(l-2),n+2k+j} \circ \mathbb{P}^x_{k,n+k+j} \circ \mathbb{P}^x_{k,n+j}  \circ   \mathbb{P}^x_{j,n} \\
 &=& \mathbb{P}^x_{k,n+kl+j} \circ \mathbb{P}^x_{k,n+k(l-1)+j} \circ \cdots \circ   \mathbb{P}^x_{k,n+2k+j} \circ \mathbb{P}^x_{k,n+j }  \circ   \mathbb{P}^x_{j,n} .
\end{eqnarray*} 
With $k:=k_0$ and $m,l,j$ such that $m \geq k_0$, $m = kl + j$ and $k_0 \leq j < 2k_0$, 
the iterated application of the contraction property in Part (3) shows that 
\[ \tilde{W}((\mathbb{P}^x_{m,n})^\ast(\,\cdot\,), (\mathbb{P}^x_{m,n})^\ast(\, \cdot\,)) 
\leq s^{ l +1 }  
\tilde{W}(\,\cdot\,,\,\cdot\,).
\]
Assertion (\ref{eq:contraction_vaserstein}) follows from this by substituting $s$ with $s^{1/2k_0}$. By a further application of Kantorovich's duality, (\ref{eq:contraction_Lipschitz}) easily follows. 
\end{proof}

As an almost immediate corollary, we obtain a family of measures. Namely, it follows from (\ref{eq:contraction_Lipschitz}) in the above Theorem, that $D_r(\mathbb{P}^x_{m,n}(g))$ tends to $0$. Moreover, $\mathbb{P}^x_{m+j,n}= \mathbb{P}^x_{m,n+j} \circ \mathbb{P}^x_{m,n} $ and $\mathbb{P}^x_{m,n+j}(\mathbf{1}) = \mathbf{1}$ imply that $\mathbb{P}^x_{m,n}(g)$ converges to a constant function. Since   
\begin{equation}
\label{eq:def_measure} \nu_x^{(n)}:\mathcal{H}_r \to \R, \; g \mapsto \lim_{m\to \infty} \mathbb{P}^x_{m,n}(g)(v)\end{equation}
is continuous by Lemma \ref{lemma_doeblin_fortet_inequality}, it follows that $\nu_x^{(n)}$ is a Borel probability measure. For ease of notation, set $\nu_x := \nu_x^{(0)}$.

\begin{theorem} \label{theo:ratio-exponential-decay} Suppose that $\varphi$ has uniformly bounded H\"older coefficients. Then there exist $C >0$ and $s \in (0,1)$ such that, for $g, h \in \mathcal{H}_r$ with $h>0$, we have for all $n \in \N$ that
\begin{equation}
\label{eq:ratio-exponential-decay} \left\| \frac{\mathcal{L}^{n}_x(g) }{\mathcal{L}^{n}_x(h)}  - \frac{\int g d\nu_x}{\int h d\nu_x} \right\|_\mathcal{H}  \leq  C s^n \left( D_r(g)  + \left|\textstyle{\frac{\int g d\nu_x}{\int h d\nu_x}}\right| D_r(h) \right) \|1/h\|_\infty.
\end{equation}
\end{theorem}

\begin{proof}
In order to show (\ref{eq:ratio-exponential-decay}), we are using  $x/y - a/b = y^{-1}((x-a) - a(y-b)/b)$ and obtain that
\begin{eqnarray*}
& & \left| \frac{\mathcal{L}^n_x(g) (w)}{\mathcal{L}^n_x(h)(w)}  - \frac{\int g d\nu_x}{\int h d\nu_x} \right|   =   
 \left|\frac{\mathbb{P}^x_{n,0}({g}) (w)}{\mathbb{P}^x_{n,0}({h})  (w)}  - \frac{\int g d\nu_x}{\int h d\nu_x}  \right| \\
 & \leq &  \frac{1}{\mathbb{P}^x_{n,0}({h})  (w)} \left(\left| \mathbb{P}^x_{n,0}({g})  (w) - {\int g d\nu_x}\right| +  \left|{\frac{\int g d\nu_x}{\int h d\nu_x}}\right|  \left| \mathbb{P}^x_{n,0}({h})  (w)  - {\int h d\nu_x} \right|\right)\\
  & \leq &  \|1/h\|_\infty \left(\left| \mathbb{P}^x_{n,0}({g})  (w) - {\int g d\nu_x}\right| +  \left|{\frac{\int g d\nu_x}{\int h d\nu_x}}\right|  \left| \mathbb{P}^x_{n,0}({h})  (w)  - {\int h d\nu_x} \right|\right). 
\end{eqnarray*}
We are now almost in position to apply (\ref{eq:contraction_vaserstein}) of Theorem \ref{label:main_prop_for_contraction}. In order to do so, first note  that $(\mathbb{P}^x_{m,0})^\ast (\nu_x^{(m)}) = \nu_x$. For $n \geq k_0$, we obtain that
\begin{eqnarray*}
\left| \mathbb{P}^x_{n,0}({g})  (w) - {\int g d\nu_x}\right| &=& \left|\int g d(\mathbb{P}^x_{n,0})^\ast(\delta_w) - \int g d(\mathbb{P}^x_{n,0})^\ast(\nu_x^{(n)}) \right| \\
&\leq & \tilde{D}(f)  s^n \tilde{W} (\delta_w,\nu_x^{(n)} ) \leq \tilde{D}(f)  s^n.
\end{eqnarray*}
Since $\tilde{D} \leq D_r$, the estimate in (\ref{eq:ratio-exponential-decay}) follows with respect to $\| \cdot \|_\infty$. The proof of the remaining assertion follows the same lines: by substituting ${\int g d\nu_x}/{\int h d\nu_x}$ with ${\mathcal{L}^n_x(g) (v)}/{\mathcal{L}^n_x(h)(v)}$ in the above estimates and applying (\ref{eq:contraction_Lipschitz}), we obtain the remaining estimate for $D_r(\mathcal{L}^n_x(g)/\mathcal{L}^n_x(h))$. 
\end{proof}

\begin{remark} Observe that the proofs of Theorems \ref{label:main_prop_for_contraction} and \ref{theo:ratio-exponential-decay} are of pathwise nature. That is, we only made use of uniform bounds for $D_r(\varphi_{f^n(x)})$ with respect to $n \geq 0$ and that $f:M \to M$ is well defined. In particular, it is not required that $f$ is invertible or measurable.  
\end{remark}

As an application of Theorem \ref{label:main_prop_for_contraction}, we now deduce unicity and continuity of the families $\{h_x\}$, $\{m_x\}$ and $\{\rho_x\}$ as defined in the introduction. As shown below, $\{m_x\}$ and $\{\rho_x\}$ are always uniquely determined by $\mathcal{L}_x^\ast(m_{f(x)})= \rho_x  m_x$, whereas unicity of  $\{h_x\}$ only holds in case of an invertible transformation. In order to deduce continuity of these families, we consider the following conditions for $f$ and $\varphi$. 

\begin{definition} \label{def:relatively_Hoelder} If $M$ is a compact topological space, we refer to $\varphi: \Sigma \times M \to \R, (w,x) \mapsto \varphi_x(w)$ as \emph{uniformly relatively $r$-H\"older continuous} if $\varphi$ is continuous and there exists $r \in (0,1)$ such that $\varphi_x \in \mathcal{H}_r$ for all $x \in M$ and $x \mapsto D_r(\varphi_x)$ is continuous. 

If, in addition, $(M,d)$ is a metric space,  we will refer to $f:M \to M$ as Lipschitz continuous if there exists $K_f>0$ such that $d(f(x),f(y))\leq K_fd(x,y)$ for all $x,y \in M$. Moreover, we will refer to $\varphi$ as \emph{$(r,\alpha)$-H\"older continuous}, if  $\varphi$ is  uniformly relatively H\"older continuous and there exists $\alpha >0 $ and $K_{\varphi}>0$ such that $|\varphi_x(w) - \varphi_y(w)| \leq K_{\varphi}d(x,y)^\alpha$ for all $x,y \in M$ and $w \in \Sigma$.    
\end{definition} 

The following proposition provides criteria for unicity, continuity and H\"older continuity of $\{m_x\}$ and $\{\rho_x\}$. 

\begin{proposition} \label{prop:continuity_1} 
Suppose that $\varphi$ has uniformly bounded H\"older coefficients.
\begin{enumerate}
\item If $\{m_x :x \in M\}$ is a family of Borel probability measures and  $\{\rho_x >0 : x \in M\}$ a family of constants such that $\mathcal{L}_x^\ast(m_{f(x)})= \rho_x  m_x$ for all $x \in M$, then these families are unique. Moreover, $m_x = \nu_x$ and $\rho_x = \nu_{f(x)}(\mathcal{L}_x(\mathbf{1}))$, with $\nu_x=\nu_x^{(0)}$ as in (\ref{eq:def_measure}). 
\item If $f$  is continuous and $\varphi$ is  uniformly relatively $r$-H\"older continuous, then  $x \mapsto \nu_x$ and $x \mapsto \rho_x$ are continuous.
\item If $f$ is Lipschitz and $\varphi$ is $(r,\alpha)$-H\"older continuous, then $x \mapsto \nu_x$ and $x \mapsto \rho_x$ are $(\min\{\alpha,\beta\})$-H\"older continuous, for  $\beta := -\log s/\log K_f$, where $s$ is equal to the contraction ratio given by Theorem \ref{label:main_prop_for_contraction}.
\end{enumerate}
\end{proposition}

\begin{proof} 
We begin with the proof of a basic estimate. For $x,y \in M$, $n \in \N$, set 
\[a_{x,y}(n) := \sup\left(\left\{\sum_{k=0}^{n-1} \left|\varphi_{f^k(x)}(\theta^k(v)) -  \varphi_{f^k(y)}(\theta^k(v))\right|   : w \in \Sigma, \theta^n(v)=w \right\}\right).\]
By the same arguments as in Lemma \ref{lemma_doeblin_fortet_inequality}, we obtain, for $g:\Sigma \to \R$, that 
\begin{eqnarray*}
|\mathbb{P}^x_{n,k}(g)(w)-\mathbb{P}^y_{n,k}(g)(w)| & \leq & 
2 (e^{a_{x,y}(n+k)}-1) \|g\|_\infty. 
\end{eqnarray*}
Hence, by Kantorovich's duality, $\tilde{W}((\mathbb{P}^x_{n,k})^\ast(\delta_w),(\mathbb{P}^y_{n,k})^\ast(\delta_w) \leq 2 (e^{a_{x,y}(n+k)}-1)$. Furthermore, if $\varphi$ is continuous, then $(\mathbb{P}^x_{n,k})^\ast(\delta_w)$ varies continuously with respect to $x$.   
If $f$ is Lipschitz and $\varphi$ is $(r,\alpha)$-H\"older, then  
\[{a_{x,y}(n)} \leq K_{\varphi}  \sum_{k=0}^{n-1} d(f^k(x), f^k(y))^\alpha \leq  \frac{K_{\varphi}(K_f^{\alpha n} -1)}{K_f^\alpha -1}d(x,y)^\alpha   \]
Hence, in this case, ${a_{x,y}(n)} \leq C K_f^{\alpha n} d(x,y)^\alpha$ for some $C\geq 1$.

\paragraph{(1)} We now prove uniqueness. Set $m_{x}^{(n)} := m_{x}\circ \theta^{-n}$  for $n \in \N$ and observe that $\mathcal{L}_x^\ast(m_{f(x)})= \rho_x  m_x$ implies that $(\mathbb{P}^x_{n,0})^\ast(m_{x}^{(n)}) = m_{x}$. 
As it was already shown above, we also have that $(\mathbb{P}^x_{n,0})^\ast(\nu_{x}^{(n)}) = \nu_{x}$. Hence, $m_x=\nu_x$ for all $x\in M$
by Theorem \ref{label:main_prop_for_contraction}. Also note that
\[ \nu_{f(x)}(\mathcal{L}_x(g)) = 
 \lim_{n \to \infty} \mathbb{P}^x_{n+1,0}(g)   \cdot \mathbb{P}^x_{n,0}(\mathcal{L}_x(\mathbf{1})) =\nu_{f(x)}(\mathcal{L}_x(\mathbf{1})) \cdot \nu_{f(x)}(g),   \] 
 implies that $\rho_x = \nu_{f(x)}(\mathcal{L}_x(\mathbf{1}))$.
\paragraph{(2 \& 3)} Let $u_n := \sup_m \{\tilde{W}((\mathbb{P}^x_{n,0})^\ast(m),(\mathbb{P}^y_{n,0})^\ast(m))\}$, where $x,y\in M$  are fixed and the supremum is taken over all Borel probability measures. By 
 Theorem \ref{label:main_prop_for_contraction}, for $n \in \N$, we have that 
\begin{eqnarray*}
&& \tilde{W}((\mathbb{P}^x_{n+k_0,0})^\ast(\nu),(\mathbb{P}^y_{n+k_0,0})^\ast(\nu)) =
\tilde{W}((\mathbb{P}^x_{n,k_0} \mathbb{P}^x_{k_0,0})^\ast(\nu),(\mathbb{P}^y_{n,k_0} \mathbb{P}^y_{k_0,0})^\ast(\nu)) \\
& \leq &  \tilde{W}((\mathbb{P}^x_{k_0,n}\mathbb{P}^x_{n,0}))^\ast (\nu),(\mathbb{P}^x_{k_0,n}\mathbb{P}^y_{n,0} )^\ast(\nu)) +  \tilde{W}((\mathbb{P}^x_{k_0,n}\mathbb{P}^y_{n,0}))^\ast (\nu),(\mathbb{P}^y_{k_0,n}\mathbb{P}^y_{n,0} )^\ast(\nu))\\
&\leq& 
u_n + s^n \tilde{W}((\mathbb{P}^x_{k_0,n})^\ast (\nu),(\mathbb{P}^y_{k_0,n})^\ast(\nu)).
\end{eqnarray*}
It follows from $\nu_x = \lim_n (\mathbb{P}^x_{n,0})^\ast(\delta_w)$ and the above estimate, that
\begin{eqnarray}
\nonumber \tilde{W}(\nu_x,\nu_y) &\leq & u_{k_0} + \sum_{n=k_0}^\infty s^n \tilde{W}((\mathbb{P}^x_{k_0,n})^\ast (\delta_w),(\mathbb{P}^y_{k_0,n})^\ast(\delta_w)) \\
\label{eq:estimate-Hoelder-measure}&\leq &  u_{k_0} + \sum_{n=k_0}^\infty s^n \min \{1,2 (e^{a_{x,y}(k_0 + n)}-1) \}
\end{eqnarray}
Since $s \in (0,1)$ and $a_{x,y}(n) \to 0$ as $x$ tends to $y$, it follows that $\nu_x$ varies continuously in $x$. If $f$ is Lipschitz and $\varphi$ is $(r,\alpha)$-H\"older, then $a_{x,y}(n) \leq C K_f^{\alpha n} d(x,y)^\alpha$. In particular, $2 (e^{a_{x,y}(k_0 + n)}-1) \geq 1$ for all $n + k_0 \geq N_0$, where 
\[ N_0 := \left\lfloor \frac{\log (\log(3/2)) - \log (C d(x,y)^\alpha)}{\log K_f^\alpha} \right\rfloor +1. \]   
Observe that $s^{N_0} = C's^{-\log d(x,y)/\log K_f} = C' d(x,y)^\beta$, for $\beta:= -\log s/\log K_f$ and some $C'>0$ and that $K_f^{\alpha N_0}d(x,y)^\alpha$ is uniformly bounded.
It then easily follows by dividing the sum in (\ref{eq:estimate-Hoelder-measure}) at $N_0$ that, for some $C''>0$,
\begin{eqnarray*}
\tilde{W}(\nu_x,\nu_y) &\leq & u_{k_0} + 
 \sum_{n=k_0}^{N_0-k_0 -1} s^n 2(e^{C K_f^{\alpha (n+k_0)} d(x,y)^\alpha}-1)    + \sum_{n=N_0 -k_0}^\infty s^n \\
 &\leq & u_{k_0} +   \frac{C}{\log 3/2}  \sum_{n=k_0}^{N_0-k_0 -1} s^nK_f^{\alpha(n + k_0)}  d(x,y)^\alpha + \frac{s^{N_0-k_0}}{1-s}\\
  &\leq& u_{k_0} +  \frac{s^{N_0}}{s^{k_0}} \left(\frac{CK_f^{\alpha N_0}d(x,y)^\alpha} {(s K_f^\alpha - 1)\log 3/2} + \frac{1}{1-s}\right)
  \leq C''(d(x,y)^\alpha + d(x,y)^\beta).
\end{eqnarray*}
This proves that $x \to \nu_x$ is H\"older continuous with index $(\min\{\alpha,\beta\})$. 
In order to prove H\"older continuity of $x \to \rho_x$, note that 
\begin{eqnarray*}
\left|\rho_x - \rho_y\right| &= & \left| \nu_{f(x)}(\mathcal{L}_x(\mathbf{1})) - \nu_{f(y)}(\mathcal{L}_y(\mathbf{1})) \right|\\
 &\leq&  \left| \nu_{f(x)}(\mathcal{L}_x(\mathbf{1})) - \nu_{f(y)}(\mathcal{L}_x(\mathbf{1})) \right| 
+ \left| \nu_{f(y)}(\mathcal{L}_x(\mathbf{1}) -\mathcal{L}_y(\mathbf{1})) \right| \\
 &\leq&  \tilde{D}(\mathcal{L}_x(\mathbf{1})) \tilde{W}(\nu_x,\nu_y) + \rho_y (e^{K_\varphi d(x,y)^\alpha} -1). 
\end{eqnarray*}
Since $\rho_y$ is uniformly bounded, $x \to \rho_x$ is $(\min\{\alpha,\beta\})$-H\"older.\end{proof}

If $f$ is invertible, then we obtain the analogue of Proposition \ref{prop:continuity_1} for  the family $\{h_x\}$. 

\begin{proposition} \label{prop:continuity_2} 
Suppose that $\varphi$ has uniformly bounded H\"older coefficients and $f$ is invertible.
\begin{enumerate}
\item If $\{h_x :x \in M\}$ is a family of strictly positive functions with uniformly bounded H\"older coefficients  and  $\{\rho_x : x \in M\}$ is a family of positive constants with $\mathcal{L}_x(h_x)= \rho_x  h_{f(x)}$ and $\int h_x d\nu_x =1$ for all $x \in M$, then these families are uniquely determined. In particular, for all $x \in M$ and $w \in \Sigma$,
\begin{equation}\label{eq:h}
 h_x(w) = \lim_{n \to \infty} \frac{\mathcal{L}^n_{f^{-n}(x)}(\mathbf{1})(w)}{\int \mathcal{L}^n_{f^{-n}(x)}(\mathbf{1}) d\nu_x} .
\end{equation}
\item If $f$ and $f^{-1}$ are continuous and $\varphi$ is uniformly relatively $r$-H\"older continuous, then  $x \mapsto h_x$ and $x \mapsto \rho_x$ are continuous.
\item If $f$ and $f^{-1}$ are Lipschitz and $\varphi$ is $(r,\alpha)$-H\"older continuous, then $x \mapsto h_x$ and $x \mapsto \rho_x$ are $(\min\{\alpha,\beta\})$-H\"older continuous (with $\beta$ as in Prop. \ref{prop:continuity_1}).
\end{enumerate}
\end{proposition}

\begin{proof} We begin with the proof of the existence of $h_x$ given by (\ref{eq:h}). Recall that there exists $C>0$ such that $1/C < \mathcal{L}^n_{f^{-n}(x)}(\mathbf{1})(v)/\mathcal{L}^n_{f^{-n}(x)}(\mathbf{1})(w) < C$ for all $v,w \in \Sigma$ and $x\in M$. Set $a_n(x):= \int \mathcal{L}^n_{f^{-n}(x)}(\mathbf{1}) d\nu_x$. It  follows from Prop. \ref{prop:continuity_1} that $a_{n+k}(x)/a_n(x) = (\rho_{f^{-n-k}(x)} \cdots \rho_{f^{-1}(x)})/(\rho_{f^{-n}(x)} \cdots \rho_{f^{-1}(x)})= a_{k}(f^{-n}(x))$.
 Hence, 
\begin{eqnarray*}
 && \left| \frac{\mathcal{L}^n_{f^{-n}(x)}(\mathbf{1})(w)}{\int \mathcal{L}^n_{f^{-n}(x)}(\mathbf{1}) d\nu_x} -   \frac{\mathcal{L}^{n+k}_{f^{-n-k}(x)}(\mathbf{1})(w)}{\int \mathcal{L}^{n+k}_{f^{-n-k}(x)}(\mathbf{1}) d\nu_x} \right| \\
& = & \frac{1}{a_n(x)}\left| \mathcal{L}^n_{f^{-n}(x)}(\mathbf{1})(w) - \mathcal{L}^{n+k}_{f^{-n-k}(x)}\left(
 {\textstyle \frac{a_n(x)}{a_{n+k}(x)}} \mathbf{1}\right)(w)  \right|
 \\
 &\leq & C \left|
 \frac{\mathcal{L}^n_{f^{-n}(x)} \left( \mathbf{1} - 
  {\textstyle \frac{1}{a_{k}(f^{-n}(x))}}
  \mathcal{L}^k_{f^{-n-k}(x)}(\mathbf{1})\right)(w)}{\mathcal{L}^n_{f^{-n}(x)}(\mathbf{1})(w)}\right|.
\end{eqnarray*}
Since $D_r(\mathcal{L}^k_{f^{-n-k}(x)}(\mathbf{1})/a_{k}(f^{-n}(x)))$ 
is uniformly bounded by Lemma \ref{lemma_doeblin_fortet_inequality}, the  sequence in (\ref{eq:h}) is a  Cauchy sequence with respect to the H\"older norm through application of  Theorem \ref{theo:ratio-exponential-decay}. Hence, $h_x$ in (\ref{eq:h}) is well defined and, in particular, $\mathcal{L}_x(h_x)= (\int \mathcal{L}_x\mathbf{1}d\nu_{x}) h_{f(x)}$. The unicity also follows from Theorem \ref{theo:ratio-exponential-decay}.  

For the proof of continuity, observe that the above  shows that the convergence in (\ref{eq:h}) is uniform with respect to $x$. The continuity of $x \mapsto h_x$ then easily follows from a $(3 \, \varepsilon)$-argument. 
The proof of  H\"older continuity follows exactly the same steps as the one in  in Proposition \ref{prop:continuity_1} and is therefore omitted.
\end{proof}

\begin{remark} 
The authors would like to point out that the novelty of this section is the method of proof based on the Vaserstein metric and the family $\{\mathbb{P}_{m,n}^x\}$. The arguments for the contraction of the 
Vaserstein metric were already known in the context of Markov operators (see \cite{Hairer}) and for normalized potentials (\cite{Stadlbauer:2015,Kloeckner:2015aa}). By considering $\{\mathbb{P}_{m,n}^x\}$, these arguments immediately can be adapted to arbitrary H\"older potentials without a priori having knowledge about the properties of $\{h_x\}$, $\{\rho_x\}$ and $\{\nu_x\}$.   
 This then allows  rather immediately to obtain results for quotients (Theorem \ref{theo:ratio-exponential-decay} above) and to recover the existence and continuity results by Denker and Gordin for the far less general situation we consider in here (see Theorems 2.6, 2.7 and 2.10 in  \cite{DenkerGordin:1999}).   
 
  However, our approach can be adapted almost in verbatim to the continuous analogues of the random bundle transformations in \cite{Stadlbauer:2015}, which are no longer skew products and are defined as follows.  Suppose that $\{(a^x_{i,j}: i,j \in \N)\}$ is a family of matrizes with entries in $\{0,1\}$ such that the family $\{\Sigma_x\}$, with 
\[\Sigma_x := \left\{(w_i): w_i \in \N,\,  
a_{w_i,w_{i+1}}^{f^i(x)} =1 \, \forall  i = 0,1,\ldots \right\} \hbox{ for } x \in M, \]  
satisfies the continuous analogue of the b.i.p. property. Furthermore, assume that $\varphi$ is a potential with uniformly bounded H\"older coefficients such that  $\|\mathcal{L}_x(\mathbf{1})\|_\infty$ is uniformly bounded. It follows from substituting the coupling in Step (2) of the proof of Theorem \ref{label:main_prop_for_contraction} by the one in \cite{Stadlbauer:2015} that  
Theorems \ref{label:main_prop_for_contraction} and \ref{theo:ratio-exponential-decay} extend to this more general situation. Furthermore, in order to extend  Propositions \ref{prop:continuity_1} and \ref{prop:continuity_2}, one additionally has to assume that $\{x \mapsto a^x_{i,j} : i,j \in \N \}$ is equicontinuous or uniformly $r$-H\"older, respectively. 
\end{remark}

\section{Stochastic and normal cocycles}\label{sec:normal} 
In this section we analyse the relation between splittings for cocycles with values in the semigroup of stochastic matrices and for the induced action, referred to as normal cocycle. For $S\colon M\rightarrow \mathcal{S}^n$, set  
\[
L_S\colon  M\times\mathbb{R}^n  \longrightarrow  M\times\mathbb{R}^n,\quad (x,v)  \longmapsto  (f(x),S_x v).
\]
Observe that the subspace generated by the unit vector $u=(1,...,1)/\sqrt{n}$ (with all entries equal to $1/\sqrt{n}$) is an invariant subspace for each $S \in \mathcal{S}^n$.  We thus obtain that the largest Lyapunov exponent $\lambda_1$ of $L_S$ is equal to zero and hence this subspace is contained in the  Oseledets subspace $E^0$.  

In order to define the induced action associated with $S \in C^0(M,\mathcal{S}^n)$, we consider two mutually orthogonal projections $P$ and $Q$ where $P$ refers to the projection onto the one-dimensional subspace $u\mathbb{R}$ and $Q$ to the projection onto the $(n-1)$-dimensional orthogonal complement $N$ of $u\mathbb{R}$. In order to define an action of $S$ on $N$, note that $S_x(N)$ is in general not contained in $N$. For this reason, we define the \emph{normal cocycle} $\hat{L}_S\colon M\times N\rightarrow M\times N$ by 
\[\hat{L}_S(x,v)=(f(x),Q \circ S_x v):=(f(x),\hat{S}_x v),\] 
or equivalently, $\hat{S}_x := Q S_x Q$, and refer to the iterates of the normal cocycle as $\hat{S}^k_x$, that is $\hat{S}^k_x := \hat{S}_{f^{k-1}(x)}\cdots\hat{S}_x$. The following elementary result now shows that the normal cocycle is a factor of the original one. 

\begin{lemma}\label{iterates2} We have  ${\hat{S}_x^k} = Q S^k_x Q $. In particular, $ S_x^k =  PS_x^{k}  + {\hat{S}_x^{k}} = P +  P S_x^{k}Q  +  {\hat{S}_x^{k}}$.
\end{lemma}
\noindent \textbf{Proof.}   The first assertion follows by induction as follows. Assume that ${\hat{S}_x^{k-1}} = Q S^{k-1}_xQ$. Then, using $P+Q= \id$, $SP=P$ and $PQ=0$, we obtain 
\begin{eqnarray*} 
Q S_x^kQ   - {\hat{S}_x^k}  
& = Q  S_{f^{k-1}(x)} (\id - Q) S_x^{k-1}Q \\
& = Q  S_{f^{k-1}(x)} P  S_x^{k-1}Q  \\
& = Q  P S_x^{k-1}Q = 0.
\end{eqnarray*}  
Hence, ${\hat{S}_x^k} = Q S^k_xQ$. Furthermore, it easily can be seen that $ P = P  S^k_{x} P  = S^k_{x} P  $. The remaining assertions then follow from
\begin{eqnarray*} 
S_x^k  = (P + Q ) S^k_{x} (P + Q)  
 =   P  + (P  + Q ) S^k_{x}  Q  =  P  + P  S^k_{x}Q    + {\hat{S}_x^k}. 
\end{eqnarray*}
~ \hfill $\boxempty$\\ 

Note that for a given $S \in C^0(M,\mathcal{S}^n)$, the normal cocycle $\hat{S}$ does not have to be in $C^0(M,\mathcal{S}^{n-1})$. However,  $\sup_{\|v\|_\infty = 1} \|Q(v)\|_\infty = 1$ implies that $\llbracket \hat{S} \rrbracket_\infty \leq 1$ and, in particular, that the Lyapunov exponents of $\hat{S}$ smaller than or equal to $0$.

By applying  Oseledets' Theorem to $\hat{L}$, we obtain the existence of $\hat{k} \in \{1,\ldots,n-1\}$, 
of exponents  $\hat{\lambda}_{1}>\ldots >\hat{\lambda}_{\hat{k}} \geq -\infty$, and
of a splitting $N_x =\hat{E}^{\lambda_1}_x \oplus \cdots \oplus \hat{E}^{\lambda_{\hat{k}}}_x$ for almost all $x \in M$ with the corresponding properties. In order to describe the relation between the splittings, we refer to $E^\lambda$ as the Oseledets subspace of exponent $\lambda$, and set 
\[E^{>\lambda}:= {\bigoplus}_{\lambda' > \lambda} E^{\lambda'}, \quad E^{\geq \lambda}:= {\bigoplus}_{\lambda' \geq \lambda}  E^{\lambda'}.\]

\begin{lemma}\label{angle} If $\lambda <0 $ is an element of the Lyapunov spectrum of $S$, then $\lambda$ is also an element of the Lyapunov spectrum of $\hat{S}$. In particular, $Q(E_x^\lambda) \subset \hat{E}_x^{\geq\lambda} \setminus  \hat{E}_x^{>\lambda}$ for almost all $x \in M$.  
\end{lemma}

\noindent \textbf{Proof.} We begin with the case of $\lambda = -\infty$. It then follows e.g. from $\|Q S_x^n(v)\| \leq \|S_x^n(v)\|$ and the definition of the Oseledets splitting that $Q(E_x^{-\infty}) \subset \hat{E}_x^{-\infty}$ for almost all $x \in M$. 
Now assume that $\lambda > - \infty$ which, in particular, implies that $S_x|_{E_x^\lambda}: {E_x^\lambda} \to {E_{f(x)}^\lambda}$ is almost surely invertible. 

For $x \in M$, let $\theta({x}):= \measuredangle (u_x, E^\lambda_x) $ refer to the non-orientated angle between $u_x$  and $E^\lambda_x$, where, with $\langle u_x, v \rangle$ referring to the standard inner product, 
\[ \measuredangle (u_x, E^\lambda_x) := \inf\left(\left\{ \cos^{-1}|\langle u_x,v \rangle|  \with v \in E^\lambda_x, \|v\|=1 \right\}\right) \in [0,\pi/2]. \]
Since $\lambda \neq 0$ it follows that $\theta({x}) >0$ on a set of full measure and hence there exists $\alpha>0$ and $A \subset M$ of positive measure with $\theta({x})>\alpha$ for all $x \in A$. As a consequence of ergodicity and invariance of $\mu$ there exist, for almost all $x \in M$, a sequence $(n_l)$ with  $n_l\nearrow \infty$  and $f^{n_l}(x) \in A$ for all $l \in \mathbb{N}$. Since $\|Q(v)\| / \|v\| \geq  \sin \alpha$ for all $y \in A$ and $v \in E^\lambda_y$, we obtain

\[
1 \geq \limsup_n \|Q S_x^{n}(v)\| / \|S_x^{n}(v)\|  \geq \sin \alpha 
\]
for almost all $x \in M$ and all $v \in  E^\lambda_x \setminus \{0\}$. When restricted to $N_x$, we have by Lemma \ref{iterates2} that $Q S_x^n = \hat{S}_x^n$. Hence, for almost all $x \in M$ and all $v \in  E^\lambda_x \setminus \{0\}$, 
\[
\limsup_{n \rightarrow{+{\infty}}}\frac{1}{n}\log{\|{S}_x^{n}(v)\|} = 
\limsup_{n \rightarrow{+{\infty}}}\frac{1}{n}\log{\|\hat{S}_x^{n}(Qv)\|}. 
\]
Hence the Lyapunov exponent of $Q(v)$ is almost surely equal to $\lambda$. In particular, $Q(v) \in \hat{E}^{\geq\lambda} \setminus  \hat{E}^{>\lambda}$.
~ \hfill $\boxempty$\\ 

\subsection{Projective hyperbolicity} The notions of dominated splittings introduced above can be extended in full analogy to normal cocycles. Furthermore, note that it is well known that domination implies continuity of the associated splitting. 
In order to relate the splittings of the cocycles, the following identity for the left hand side of (\ref{ds}) will be essential.
\begin{eqnarray}
 \frac{\|S^{m}_x|_{F^2_{x}}\|}{\mathfrak{m}(S^{m}_x|_{F^1_{x}})} 
 & = & \sup_{w \in F^2_{x}\setminus\{0\}}  \frac{\|S^{m}_x(w)\|}{\|w\|} 
    {\sup_{ v \in S^{m}_x(F^1_{x}\setminus\{0\}) }  \frac{\|(S^{m}_x|_{F^1_{x}})^{-1} (v)\|}{\|v\|} }  \label{sup-inf}  \\ 
& = & \sup_{w \in F^2_{x}\setminus\{0\}}  \frac{\|S^{m}_x(w)\|}{\|w\|} 
{\sup_{ v \in F^1_{x}\setminus\{0\} }  \frac{\|v\|}{\|S^{m}_x(v)\|} } \nonumber \\
& = & \frac{\sup_{w \in F^2_{x}\setminus\{0\}}  \frac{\|S^{m}_x(w)\|}{\|w\|}}{\inf_{ v \in F^1_{x} \setminus\{0\}}  \frac{\|S^{m}_x(v)\|}{\|v\|} } . \nonumber
\end{eqnarray}
Hence, if the splitting is dominated, then the above is smaller than $1/2$ and, in particular, the minimal Lyapunov exponent in $F^1$ is bounded away from the maximal exponent of $F^2$.

We now relate these two concepts of dominated splitting in the following lemma, which is formulated for the slightly more general case of a compact and $f$-invariant subset  $\Lambda$ of $M$. In here, we refer to an invariant bundle $V \subset M \times \R^n$ as an \emph{eventually contracting} bundle if there  exists $k \in \N$ with $\sup_{x \in M} \|S^k_x|_{V_x}\| <1$. 
The hypothesis of $\Lambda$ being a compact is crucial in order to apply some abstract invariant manifold paraphernalia as in \cite{HPS}. 

\begin{lemma}\label{equivalence}
Let $\Lambda \subset M$ be compact and $f$-invariant. Then the cocycle $S$ has a dominated splitting  over $\Lambda$  if the normal cocycle $\hat{S}$ admits a dominated splitting over $\Lambda$. Furthermore, if $\hat{S}$ is eventually contracting, then $\dim(E^0_x)=1$ for all $x \in \Lambda$ and $E^0_\Lambda \oplus E^{<0}_\Lambda$ is a dominated splitting.
\end{lemma} 
\noindent \textbf{Proof.} We begin with the proof of the first statement. So assume that $N^1\oplus N^2$ is a given $\hat{m}$-dominated splitting for $\hat{S}$, that is  ${\|\hat{S}^{\hat{m}}_x|_{N^2_{x}}\|}/{\mathfrak{m}(S^{\hat{m}}_x|_{N^1_{x}})}\leq \frac{1}{2}$. Set $F^1 := u\mathbb{R}_{\Lambda}\oplus N^1$ and note that the sequence
\[
\begin{CD} 
0 	@>>> u\R_{\Lambda}   @>>>  u\R_{\Lambda}  \oplus N^2  @>Q>> N^2  @>>> 0
\end{CD}
\]
is exact. Moreover, since $\hat{S}$ is a factor of $S$, the following diagram commutes.
\[
\begin{CD} 
0 	@>>> u\R_{\Lambda}   @>>>  u\R_{\Lambda}  \oplus N^2  @>Q>> N^2  @>>> 0 \\
@.       @VV{S|_{u\R}}V           @VV{S|_{ u\R \oplus N^2 }}V     @VV\hat{S}V  @. \\
0 	@>>> u\R_{\Lambda}   @>>>  u\R_{\Lambda}  \oplus N^2  @>Q>> N^2  @>>> 0
\end{CD}
\]

We now construct an $S$-invariant subbundle $F^2\subset N^2 \oplus u\mathbb{R}_{\Lambda}$ following ideas in the proof of Lemma 2.18 in \cite{HPS}. In here, it is sufficient to construct a continuous family of linear maps $\sigma_x: N^2_x \to u\R_x$ such that $S_x(\sigma_x(v),v) = (\sigma_{f(x)}(\hat{S}_x(v)),\hat{S}_x(v))$ for all $x \in M$, or equivalently, by Lemma \ref{iterates2},
\[  \iota_{x} \sigma_x + PS_x  = \sigma_{f(x)}\circ \hat{S}_x,\]
where $\iota_{x}^k$ stands for the identity seen as an isomorphism from $\R^n_x$ to $\R^n_{f^{k}(x)}$.
We now construct $\sigma$ as a fixed point of the map defined by 
\begin{equation} \label{eq:contraction-Banach}
\sigma_x \mapsto \iota^{-\hat{m}}_{f^{\hat{m}}(x)} \left(\sigma_{f^{\hat{m}}(x)}\circ \hat{S}^{\hat{m}}_x - P S^{\hat{m}}_x\right)
\end{equation}
on the bundle of linear maps  $L_\Lambda(N^2_\Lambda,u\R_\Lambda)$ equipped with the supremum norm. 
It follows from equation $(\ref{sup-inf})$ and the property that $N^1\oplus N^2$ is $\hat{m}$-dominated splitting, that $\llbracket \hat{S}^{\hat{m}}|_{N_\Lambda^2}\rrbracket \leq 1/2$.
Therefore, the map defined in (\ref{eq:contraction-Banach}) is a contraction and an adaption of the Banach fixed point theorem  (see \cite[Theorem 3.1]{HPS}) shows that the above map has a fixed point $\sigma$ in $L_\Lambda(N^2_\Lambda,u\R_\Lambda)$. This proves the existence of $\sigma$ where the corresponding $S$-invariant bundle is defined by $F_x^2 := \{(\sigma_x(v),v) \with v \in N_x^2\}$.

It remains to show that $F^1 \oplus F^2$ is a dominated splitting. Note that it follows from domination with respect to $\hat{S}$ that $\inf_{x \in \Lambda} \measuredangle(N^1_x, N^2_x)>0$. Since  $v \mapsto (\sigma_x(v),v)$ is continuous and $\Lambda$ is compact, we hence have that  $\measuredangle(u_x, F_x^2)$ is uniformly bounded away from $0$. In particular, for some norm  $\|\cdot\|_\ast$ equivalent to the Euclidean norm, there exists $C>1$ with $C^{-1} \|v\| \leq \|Q(v)\| \leq \|v\|$ for all $v \in F_x^2$ and $x \in \Lambda$ where the involved constant only depends on the minimal angle and the norm $\|\cdot \|_\ast$. Hence, for $k \in \N$ and with $\mathfrak{m}_\ast$ referring to the co-norm induced by $\|\cdot \|_\ast$,
\begin{eqnarray*}
\frac{\|{S}^{k}_x|_{F^2_{x}}\|_\ast}{\mathfrak{m}_\ast({S}^{k}_x|_{F^1_{x}})} 
&= \sup_{w \in F^2_{x}}  \frac{\|{S}^{k}_x(w)\|_\ast}{\|w\|_\ast} 
{\sup_{ v \in  F^1_{x} }  \frac{\| v\|_\ast}{\|{S}^{k}_x(v)\|_\ast} } \\
&\leq C \sup_{w \in F^2_{x}}  \frac{\| Q S^{k}_x(w)\|_\ast}{\|Q w\|_\ast} 
{\sup_{ v \in N^1_{x}, t \in \R }  \frac{\|v + t u_x\|_\ast}{\| t u_x + S^{k}_x(v)\|_\ast} }. 
\end{eqnarray*}
Since the property of admitting a dominated splitting is independent of the norm on $\R^n = u \R \oplus N$, we may choose $\|v + t u_x\|_\ast :=  \|v\| + |t| $. 
For fixed $v \in N^1_x$ it is then easy to verify that the second term in the above estimate achieves its global maximum for $t = - \|PS_x^k v \|$. We hence obtain 
\[  
{\sup_{ v \in N^1_{x}, t \in \R }  \frac{\|v + t u_x\|_\ast}{\| t u_x + S^{k}_x(v)\|_\ast} }
\leq \sup_{ v \in N^1_{x}} \frac{\|v\| + \|PS_x^k v \|}{\|QS^{k}_x(v)\|}  
\leq 2 \sup_{ v \in N^1_{x} }  \frac{\|v\|}{\|\hat{S}^{k}_x(v)\|}. 
\] 
By combining the two estimates above, we arrive at
\[{\|{S}^{k}_x|_{F^2_{x}}\|_\ast}/{\mathfrak{m}_\ast({S}^{k}_x|_{F^1_{x}})} \leq 2C  {\|\hat{S}^{k}_x|_{N^2_{x}}\|_\ast}/{\mathfrak{m}_\ast(\hat{S}^{k}_x|_{F^2_{x}})} \]
and hence $F^1 \oplus F^2$ is an $m$-dominated splitting for $S$ for each multiple $m \in \N$ of $\hat{m}$ with $2C/2^{m/\hat{m}} \leq 1/2$.

For the proof of the second assertion, assume that $S \in \mathcal{S}$ satisfies $\sup_{x \in \Lambda} \|\hat{S}^k_x\| =: \rho <1$. By the same arguments, it  follows from $\rho<1$ that 
\[ L_\Lambda(N_\Lambda,u\R_\Lambda) \to L_\Lambda(N_\Lambda,u\R_\Lambda), \quad \sigma_x \mapsto \iota^{-1}_{f(x)} \left(\sigma_{f(x)}\circ \hat{S}^k_x - P S^k_x\right) \]
is a contraction and hence has a unique fixed point. With $\sigma$ referring to the fixed point and $F^2_x:=\{(\sigma_x(v),v) \with v \in N\}$, it then follows as above from continuity of $\sigma$ and compactness of $\Lambda$ that the angle of $u_x$ and $F^2_x$ is uniformly bounded  away from $0$ for all $x \in \Lambda$. This then implies that the Lyapunov exponents of $F^2_\Lambda$ are bounded by $\log \rho/k$, and, in particular, that  $E^{<0}_\Lambda = F^2_\Lambda$ and $E^0_\Lambda = u\R_\Lambda$. Finally, observe that the arguments above for proving that $F^1\oplus F^2$ is dominated apply in verbatim to $E^0_\Lambda \oplus E^{<0}_\Lambda$ (and even can be simplified using $\mathfrak{m}_\ast({S}_x|_{E^0_{x}})=1$). 
~ \hfill $\boxempty$\\ 

As a corollary of the construction of $F^2$ in the above result, we obtain that $E^0 \oplus E^{<0}$ is a dominated splitting for a generic cocycle.

\begin{corollary}\label{corollary:almost_rpf} If $M$ is compact, then there exists an open and dense subset $\mathcal{R}$ of $\mathcal{S}$ such that $\dim(E^0)=1$, $E^{<0}$ is an eventually contracting bundle and $E^0 \oplus E^{<0}$ is a dominated splitting for each $S \in \mathcal{R}$.
\end{corollary}

\noindent\textbf{Proof.}  Observe that the map $S \mapsto \llbracket \hat{S}^k\rrbracket$ 
is continuous for all $k \in \N$. Hence, 
\[ \mathcal{R} := 
\bigcup_{k=1}^\infty \{S \in \mathcal{S} \with \llbracket \hat{S}^k\rrbracket  <1 \}  = \{S \in \mathcal{S} \with \exists k \in \N \hbox{ s.t. }\llbracket \hat{S}^k\rrbracket  <1 \}  \]
is open. For the proof of the corollary, it therefore remains to show that $\mathcal{R}$ is dense and then apply Lemma \ref{equivalence} to obtain that $\dim(E^0)=1$ and $E^0 \oplus E^{<0}$ is dominated. 

In order to prove that $\mathcal{R}$ is dense, we show that
\[\Psi: \mathcal{S} \times [0,1] \to \mathcal{S}, \; (S,\rho) \mapsto S^{(\rho)} := P + \rho(PSQ + \hat{S}).
\] 
is an isotopy with $S^{(\rho)} \in \mathcal{R}$ for all $\rho <1$ which continuously connects $S^{(1)} = S\in \mathcal{S}$ with $ S^{(0)} = P$. We begin with verifying that $S^{(\rho)} \in \mathcal{S}$. Note that, for $\rho \in [0,1]$, 
\begin{equation} \label{eq:u-invarince}S^{(\rho)}(u) = (P + \rho(PSQ + \hat{S}))u =u,\end{equation}
and, with $\R^n_+$ referring to $\{(v_1, \ldots, v_n)\in \R^n: v_i \geq 0,\;  i=1, \ldots, n\}$,
\begin{eqnarray}\nonumber
 v \in \R^n_+ & \Rightarrow {S}(v)  \geq 0  \Rightarrow P(v) \geq - (PSQ + \hat{S})(v) \\
 \label{eq:invariance_of_the_positive_cone}
 & \Rightarrow P(v) \geq - \rho(PSQ + \hat{S})(v)
  \Rightarrow \hat{S^{(\rho)}}(v) \in \R^n_+. 
\end{eqnarray}
It then follows from (\ref{eq:u-invarince}) and (\ref{eq:invariance_of_the_positive_cone}) that  $S^{(\rho)} \in \mathcal{S}$. In order to see that $\Psi$ is continuous, note that, for $S,T \in \mathcal{S}$ and $\rho,\rho' \in [0,1]$, 
\begin{eqnarray*} 
\llbracket S^{(\rho)} - T^{(\rho')}\rrbracket_2  & = \llbracket  S^{(\rho)} - S^{(\rho')} +S^{(\rho')} - T^{(\rho')} \rrbracket_2  \\
& \leq |\rho - \rho'| \llbracket  SQ \rrbracket_2  + \rho' \llbracket (S-T)Q\rrbracket_2  \leq |\rho - \rho'| + \rho' \llbracket S-T\rrbracket_2 . 
\end{eqnarray*}
Finally, Lemma \ref{iterates2} implies that $\hat{S^{(\rho)}} = \rho \hat{S}$. Hence, $\llbracket\hat{S}^k \rrbracket_\infty \leq \rho^k$ and, in particular,  $S^{(\rho)} \in \mathcal{R}$ for all $\rho <1$. 
~ \hfill $\boxempty$\\ 

We now give a short discussion of the continuity of the pertubation in the above corollary from the viewpoint of the Oseledets splitting in case $\dim(E^0)=1$.
Let $E^\lambda_x(\rho)$ refer to the Oseledets subspace of exponent $\lambda$ with respect to $S^{(\rho)}$, and $\hat{E}^\lambda_x(\rho)$ to the one of $\hat{S^{(\rho)}}$, respectively. 
\begin{proposition} Assume that $\dim(E^0)=1$. Then, for almost every $x \in M$ and each negative Lyapunov exponent $\lambda$, the restriction $\Phi^\rho|_{ E_x^\lambda}$ is an isomorphism from $E_x^\lambda$ to $E_x^{\lambda+ \log \rho}(\rho)$. Furthermore, 
$\lim_{\rho \to 1 - } \Phi^\rho(v) = v$ for each $v \in E_x^\lambda$. 
\end{proposition}
\noindent \textbf{Proof.}   Throughout this proof, we assume that $\lambda <0$ is a negative Lyapunov exponent of $S$.
It then follows from $\dim(E^0)=1$ and Lemma \ref{angle}, that $Q|_{E_x^\lambda}$ is an isomorphism onto $\hat{E}_x^\lambda$. Since $\hat{E}_x^\lambda = \hat{E}^{\lambda+\log \rho}_x(\rho)$, the first assertion follows from the properties of $\sigma^{(\rho)}$. It hence remains to show that $\lim_{\rho \to 1 - } \Phi^\rho(v) = v$, or equivalently, $\lim_{\rho \to 1 - } \sigma^{(\rho)}(Qv) = Pv$, for $v \in E_x^\lambda$. Therefore, a straightforward induction argument combined with the identity $S^{(\rho)}=P + \rho SQ$ gives that 
\begin{equation}\nonumber
\label{eq:iterates_of_the_deformation}
(S^{(\rho)})^n = \sum_{k=0}^n \rho^k \iota^{n-k}P(SQ)^k + \rho^n \hat{S}^n, \quad \forall n \in \N. 
\end{equation}
For $\rho =1$ and $ v \in {E}^{\lambda}_x$ we therefore obtain that 
\[\lim_{n\to \infty} S^n(v) = \lim_{n\to \infty} (S^{(\rho)})^n(v) = 0\]
 and  $\lim_n \hat{S}^n v =0$. In particular, by identifying the range of $\iota^k P$ with $\R$, we obtain     
\begin{equation}\label{eq:iterates_of_the_deformation_I}
P(v) = -  \sum_{k=1}^\infty P(SQ)^k(v).
\end{equation}
For $\rho < 1$, it follows from the above that $\Phi^\rho(v)=(\sigma^{(\rho)}(Qv),Qv) \in {E}^{\lambda+\log \rho}_x(\rho)$. Therefore, $\lim_n (S^{(\rho)})^n(\Phi^\rho(v)) = 0$ and
\begin{equation}\label{eq:iterates_of_the_deformation_II}
\sigma^{(\rho)}(Qv) = -  \sum_{k=1}^\infty \left( P(SQ)^k(v)\right) \rho^k .
\end{equation}
Set  $a_k:= P(SQ)^k(v)$ for $k \geq 1$. It follows from (\ref{eq:iterates_of_the_deformation_I}) that $\sum_{k \geq 1} a_k$ converges. Hence, by Abel's continuity theorem, 
$\lim_{\rho \to 1-} \sum_{k\geq 1} a_k \rho^k = \sum_{k\geq 1} a_k$.
The remaining assertion hence follows from (\ref{eq:iterates_of_the_deformation_II}).
~ \hfill $\boxempty$\\ 

\section{Accessibility and the proof of Theorem C}\label{proofs}
The proof of our main result for stochastic cocycles is based on the dichotomy established in \cite{BV} for cocycles with values in an accessible group. Therefore, we  recall the definition of accessibility given in \cite[Definition 1.2]{BV}.

\begin{definition} \label{def:accessible} 
Assume that $\mathfrak{G}$ is an embedded submanifold of $\mathrm{GL}_n(\R)$ with or without boundary. We say that $\mathfrak{G}$ is accessible if for all $C > 0$ and $\epsilon > 0$, there are $m \in\mathbb{N}$ and $\alpha>0$ satisfying the following properties. Given $\xi$ and $\eta$ in
the projective space $\mathbb{R}\emph{\textbf{P}}^{n-1}$ with $\measuredangle(\xi,\eta)<\alpha$ and $S_0,\ldots, S_{m-1}$ in $\mathfrak{G}$ with $\|{S_i}{}^{\pm 1}\|<C$, there exist $R_0,\ldots, R_{m-1}$ in $\mathfrak{G}$ such that $\|R_i-S_i\|<\epsilon$ and
\[{R}_{m-1}\circ ...\circ {R}_0(\eta)={S}_{m-1}\circ ...\circ {S}_0(\xi).\]
\end{definition}

The idea of proof of the main theorem is to employ the dichotomy obtained by Bochi and Viana in \cite{BV} (see Lemma \ref{corollaryBV} below) to cocycles with values in $Q\mathcal{S}^\ast Q$, where  $\mathcal{S}^\ast$ refers to the semigroup of invertible elements in $\mathcal{S}$.
 For this purpose, it is necessary to verify that $\mathcal{S}^\ast$ is dense in $\mathcal{S}$ and that $Q\mathcal{S}^\ast Q$ is accessible, that is $Q\mathcal{S}^\ast Q$ is an embedded submanifold with boundary in $\mathrm{GL}_n(\R)$ and that the approximation property in $\mathbb{R}\emph{\textbf{P}}^{n-1}$  as stated above holds. The final result then follows applying the results obtained in the previous section.

\begin{remark} \label{remark:non-accessible} The result of Bochi and Viana can not be applied directly, since  $\mathcal{S}^\ast$ is not accessible. This can be seen by the following argument. Assume that ${S}_0,\ldots ,{S}_{k}$ are elements in $\mathcal{S}^\ast$. Since $u$ is an eigenvector for all elements of $\mathcal{S}^\ast$, it  follows from invertibility that  ${S}_0 \cdots {S}_{k}(v) =u $ if and only if $v=u$. Hence, ${R}_0 \cdots {R}_{k}(\xi) \neq u $ for all ${R}_0, \ldots, {R}_{k} \in \mathcal{S}^\ast$ and $\xi\not=u$ which gives that $\mathcal{S}^\ast$ is not accessible.  
\end{remark}

We now proceed with the proof. In  order to do so, observe that $\mathcal{S}$ can be identified with the $n$-fold product of the unit simplex in $\R^n$, that is, with $\Delta_n:=\{(x_1,\ldots,x_n) \in [0,1]^n \with x_1+\ldots+ x_n=1\}$, we have  
\[ \mathcal{S} = \{ (a_{ij})_{1\leq i,j \leq n} \with (a_{i1},a_{i2},\ldots, a_{in})\in \Delta_n \;\forall i=1,\ldots ,n\} \cong (\Delta_n)^n.\]
In particular,  $\mathcal{S}$ is an embedded submanifold with boundary of the manifold  of real valued matrices $\mathrm{M}_{n\times n}(\R)$ of dimension $n\times n$ with respect to the canonical manifold structure.

\begin{lemma} \label{lem:open_dense}
$\mathcal{S}^\ast$ is dense and open in $\mathcal{S}$.
\end{lemma}
\noindent \textbf{Proof.} First observe that a standard argument using the continuity of the determinant implies that  $\mathcal{S}^\ast$ open in $\mathcal{S}$. In order to prove the density of $\mathcal{S}^\ast$, fix  $S \in \mathcal{S} \setminus \mathcal{S}^\ast$ with $v^{(1)},\ldots, v^{(n)} \in \Delta_n$ referring to the row vectors of $S$. 
Furthermore, set $k := \mathrm{Rank}(S) =\dim(\mathrm{Span}(\{v^{(1)},\ldots, v^{(n)} \}))$ and note that non-invertibility implies that $k<n$. 
Hence, there exists $l$ with 
\[v^{(l)} \in V, \; V := \mathrm{Span}(\{v^{(1)},\ldots,v^{(l-1)},v^{(l+1)},\ldots , v^{(n)} \}).\] 
Since $k<n$, it follows that for each $\epsilon > 0$, there exists $w \in \Delta_n \setminus V$ with $\|v^{(l)} - w\| < \epsilon$. By substituting $v^{(l)}$ with $w$, one obtains $T \in \mathcal{S}$ close to $S$ with $\mathrm{Rank}(T)=\mathrm{Rank}(S)+1$.
The remaining assertion then follows by induction until one obtains an element of $\mathcal{S}$ of full rank.
~ \hfill $\boxempty$\\ 

In order to prove that $Q\mathcal{S}^\ast Q$ is an embedded submanifold with boundary, we consider the following affine subspaces of $\mathrm{M}_{n \times n}(\R)$ and $\R^n$. 
\begin{eqnarray*}
\hat{\mathcal{S}} & := &   \left\{ (a_{ij}) \in \mathrm{M}_{n \times n}(\R) \with {\textstyle \sum_{j=1}^n a_{ij}= 1} \;\forall i=1, \ldots, n \right\},  \\
\mathcal{A} & := & \left\{ (a_{ij}) 
\in \mathrm{M}_{n \times n}(\R)
 \with {\textstyle \sum_{j=1}^n a_{i^\prime j}= \sum_{i=1}^n a_{i j^\prime}= 0 } \;\forall i^\prime, j^\prime =1, \ldots, n \right\}, \\
\hat{\Delta}_n  &:=&  \left\{ (v_{j}) \in \R^n \with {\textstyle \sum_{j=1}^n v_{j}= 1} \right\}.
\end{eqnarray*}
Moreover,  consider the homomorphisms of vector spaces
\begin{eqnarray}
\nonumber \kappa &:&  \mathrm{M}_{n\times n}(\R) \to \R^n, \; (a_{ij})\mapsto \left({\textstyle \frac{1}{n} \sum_{i=1}^n a_{i1}, \ldots, \frac{1}{n}  \sum_{i=1}^n a_{in} }\right),\\
\label{eq:Theta} \Theta &:& \hat{\mathcal{S}}  \to \mathcal{A}  \times \hat{\Delta}_n, \; S \mapsto (QS,\kappa(S)). 
\end{eqnarray}
It is easy to see, for $S \in \hat{\mathcal{S}}$ and with $(\cdot)_{ij}$ referring to the coordinate $(i,j)$, that $(QS)_{ij}=a_{ij} - (\kappa(S))_j$, which then implies that $(\Theta^{-1}(A,v))_{ij} = (A)_{ij} + v_j$ and that $\Theta$ is an isomorphism of vector spaces. 

\begin{lemma} \label{lem:embedded_manifold}
$Q\mathcal{S}^\ast Q$ is an embedded submanifold with boundary of $\mathrm{GL}(N)$.
\end{lemma}

\noindent \textbf{Proof.} We begin identifying $Q\mathcal{S}$ with $\mathcal{M}$, where 
\[ \mathcal{M} := \left\{ (a_{ij}) \in \mathcal{A} \with {\textstyle \sum_{j=1}^n \min_{i=1,\ldots,n}} a_{ij} \geq -1 \right\}.\]
So assume that $S=(s_{ij}) \in \mathcal{S}$. It then follows from the definition of $\Theta$ that $QS \in \mathcal{A}$. Furthermore, since $S$ is a stochastic matrix, we have that $1 \geq \sum_j \min_{i}a_{ij} \geq 0$. It then follows from $(QS)_{ij}=a_{ij} - (\kappa(S))_j$ that $\sum_j \min_{i}(QS)_{ij} \geq -1$.
Hence, $QS \in  \mathcal{M}$ and, in particular, $Q\mathcal{S}  \subset \mathcal{M}$. 
Now assume that $(b_{ij})\in \mathcal{M}$ and choose $v \in \Delta_n$ with $-(v)_j \leq \min_j b_{ij}$ for all $j =1, \ldots, n$. It then follows that $(\Theta^{-1}((b_{ij}),v))_{ij} = b_{ij} + v_j \geq 0$. Furthermore, since $\sum_j b_{ij} + v_j =1$, we also have that $b_{ij} + v_j \leq 1$. This then implies that $\Theta^{-1}((b_{ij}),v) \in \mathcal{S}$ and, in particular, that $Q\mathcal{S}  = \mathcal{M}$.      

In order to see that $\mathcal{M}$ is a submanifold, note that 
$\mathcal{M}$ can be written as 
\[\mathcal{M} = \mathcal{A} \cap \bigcap_{\sigma \in \mathfrak{S}_n} \left\{ (b_{ij}) \in \mathrm{M}_{n \times n}(\R) \with  {\textstyle \sum_{j=1}^n b_{\sigma(j) j} \geq -1} \right\}, \]   
where $\mathfrak{S}_n$ stands for the symmetric group on $\{1, \ldots,n\}$. Hence, $\mathcal{M}$ is equal to the intersection of the linear subspace $\mathcal{A}$ with finitely many affine half spaces, which proves that $\mathcal{M}$ is an embedded submanifold of $\mathrm{M}_{n \times n}(\R)$. 
The assertion then follows from $QSQ =QS$ for all $S \in \mathcal{S}$, the fact that $\mathcal{A}$ can be identified with the endomorphisms of $N$ and $QS \in \mathrm{GL}(N)$ for all $S \in \mathcal{S}^\ast$.
~ \hfill $\boxempty$\\

We now proceed with the proof that $Q\mathcal{S}^\ast Q$ is accessible. For this purpose, we analyse the orbit of an element of the positive cone $\R_+^n:= \{ (v_1,\ldots,v_n) \in \R^n \with v_i\geq 0 \; \forall i  \}$ under the following subsets of $\mathcal{S}$ which are defined by, for given $\epsilon>0$ and with $\delta_{ij}$ referring to Kronecker's $\delta$-function,
\[\mathcal{S}_\epsilon := \left\{ (s_{ij})\in \mathcal{S} \with \max_{i,j}|s_{ij} - \delta_{ij}| \leq \epsilon \right\}.\] 
 
Note that, for $\epsilon >0$ sufficiently small, each element of $\mathcal{S}_\epsilon$ is invertible.  Moreover, for $v = (v_i) \in  \R^n$, set $\alpha(v) := \min_{j} v_j$ and $\beta(v) = \max_{j} v_j$. Using the elementary fact that 
$ [\alpha,\beta] = \{ t \alpha + (1-t) \beta \with t \in [0,1]\} $ 
we then obtain a precise description of the orbit 
\[ \mathcal{S}_\epsilon(v) = \left\{ S(v)  \in   \R^n \with  S \in  \mathcal{S}_\epsilon \right\} \]
of an element $v \in \R^n$ under $\mathcal{S}_\epsilon$. That is, $w = (w_i) \in  \mathcal{S}_\epsilon(v)$ if and only if, for all $i \in \{1, \ldots,n\}$,
\begin{equation} \label{eq:orbit_of_v}
 \epsilon(\alpha(v) - v_i) \leq w_i - v_i \leq \epsilon (\beta(v)- v_i).
 \end{equation}
For the proof of the main lemma, recall that $\|v\|_\infty \leq \| v \| \leq \sqrt{n}  \| v\|_\infty $, for each $v \in \R^n$, and set
\begin{equation} \label{eq:constant-projective}
C := (n-1)\sqrt{n} + 2n(n-1)^2
\end{equation}

\begin{lemma} \label{lemma:projective_equality}
For all $\epsilon>0$, ${v},{w} \in N=Q(\R^n)$ with $\|{v}\| = \|{w}\| =1$ and $\|{v} - {w}\|\leq \epsilon$ there exist $S \in \mathcal{S}_{C\sqrt{\epsilon}}$ and $t >0$ such that ${S}{v} = t {w}$.   
\end{lemma}

\noindent \textbf{Proof.}   The first step in here is to obtain lower bounds for $\beta({v})$ and $|\alpha({v})|$ using the fact that $v,w \in N$. Since the equivalence of norms implies that $ \|{v}\|_\infty \geq \|{v}\|/\sqrt{n}$, either $\beta({v})$ or $|\alpha({v})|$ is bigger than or equal to ${1}/{\sqrt{n}}$. Furthermore, it follows from  ${v} = ({v}_1, \ldots,  {v}_n) \in N$ that $P(v)=(\sum_i {v}_i)/\sqrt{n} = 0$. We hence obtain, e.g. by analysing the worst case, that 
 \[ \min \{\beta({v}), |\alpha({v})| \}\geq \frac{1}{(n-1)\sqrt{n}} =: {D}^{-1}.  \]
The second step is to shrink $w$ in a controlled way by some $t \in (0,1]$ such that $t w \in \mathcal{S}_{\delta_2}(v)$ for some $\delta_2>0$. 
In order to determine $t$, note that $\|{v} - {w}\|_\infty \leq \|{v} - {w}\|\leq \epsilon$ implies that $|\alpha(v) - \alpha(w)| \leq \epsilon $ and  
$|\beta(v) - \beta(w)| \leq \epsilon $. For $\delta_1>0$ to be specified later, let
\begin{eqnarray*} t := & \frac{\textstyle \max(\{\beta({v}), |\alpha({v})|\}) - \delta_1}{\textstyle \max(\{\beta({v}), |\alpha({v})|\}) + \epsilon}\\ 
= &1- \frac{\textstyle \delta_1 + \epsilon}{\textstyle \max(\{\beta({v}), |\alpha({v})|\}) + \epsilon} \geq 1- D(\delta_1 + \epsilon). \end{eqnarray*}
For this choice of $t$, we have $\alpha(v) + \delta_1 \leq \alpha(t w)$, $\beta(t w)  \leq \beta(v) - \delta_1$ and
\[ \|v - t w\|_\infty  \leq \epsilon + (1-t)\|w \|_\infty \leq  \epsilon + D(\delta_1 + \epsilon).  \]
For $\delta_1:= \sqrt{\epsilon}$ and  $\delta_2:= (D + 2D^2)\sqrt{\epsilon}$ we then obtain  $\epsilon \leq \delta_1 \delta_2 $ and 
\begin{eqnarray*} 
 \epsilon + D(\delta_1 + \epsilon) \leq (1+2D)\delta_1 = \delta_2/D.
\end{eqnarray*}
In order to verify that  $t w \in \mathcal{S}_{\delta_2}(v)$ we consider the following cases for $i=1, \ldots,n$. 
\begin{enumerate}
  \item 
  If $v_i \geq t w_i$ and $v_i \geq 0$, or  $v_i \leq t w_i$ and $v_i \leq 0$, then 
\[|v_i - t w_i| \leq \delta_2/D \leq 
 \left\{\begin{array}{ll}
 \delta_2(v_i - \alpha(v)) &: v_i \geq 0,\\
 \delta_2(\beta(v)-v_i) &: v_i \geq 0. 
 \end{array}\right.\]
  \item If $0 \leq v_i \leq t w_i$, then $v_i \leq \beta(v) - \delta_1$, and if $0 \geq v_i \geq t w_i$, then $v_i \geq \alpha(v) + \delta_1$. In particular,  
  \[|v_i - t w_i|  \leq \epsilon \leq \delta_1 \delta_2  \leq     
 \left\{\begin{array}{ll}
 \delta_2(\beta(v)-v_i) &: 0 \leq v_i \leq t w_i,\\
 \delta_2(v_i - \alpha(v)) &: 0 \geq v_i \geq t w_i. 
 \end{array}\right.\]
\end{enumerate}
It now follows from (\ref{eq:orbit_of_v}) that $t w \in \mathcal{S}_{\delta_2}(v)$.
~ \hfill $\boxempty$\\ 

As an immediate consequence, we obtain the following statement about the accessibility of $\mathcal{S}$.
\begin{lemma}\label{main} For all $\epsilon > 0$, $S \in\mathcal{S}$ and $x,y \in N\setminus \{0\}$ with $\measuredangle(x,y) \leq \epsilon^2/(C^2)$ and $C$ defined in (\ref{eq:constant-projective}), there exists $R\in \mathcal{S}$ and $\lambda \in \R$ such that $Sx=\lambda Ry$ and 
$ \llbracket S-R \rrbracket \leq  \epsilon$.
In particular, $Q\mathcal{S}^\ast Q$ is accessible.
\end{lemma}

\noindent \textbf{Proof.}   Assume without loss of generality that $\|x\|=\|y\|=1$. Since $\|x-y\| \leq\measuredangle(x,y)$, we have  $\|x-y\| \leq \epsilon^2/(C^2)$. It hence follows from the above Lemma that there exist $T \in \mathcal{S}_{\epsilon}$ and $t \in (0,1]$ such that $Tx = t y$. Now set $R := S \circ T$ and observe that 
\[ \llbracket S-R \rrbracket = \llbracket S(\id - T)\rrbracket \leq \llbracket S \rrbracket \llbracket \id - T \rrbracket \leq   \epsilon.  \]
In order to prove accessibility, note that $\mathcal{S}_{\epsilon} \subset \mathcal{S}^\ast $ for $\epsilon >0$ sufficiently small. Hence, if $S \in \mathcal{S}^\ast $, then $R = S \circ T \in \mathcal{S}^\ast$ 
for $\epsilon$ sufficiently small. 
It now follows from Lemma \ref{lem:embedded_manifold} and  $\mathrm{GL}(N) \cong \mathrm{GL}_{n-1}(\R)$ that $Q\mathcal{S}^\ast Q$ is accessible.
~ \hfill $\boxempty$\\ 

As the final ingredient of the proof of theorem C, recall the Bochi-Viana dichotomy established in \cite{BV} for accessible cocycles.  
\begin{lemma}[Corollary 1 in \cite{BV}]  \label{corollaryBV}
Assume that $(f,\mu)$ is ergodic and $\mathfrak{G} \subset GL_n(\mathbb{R})$ is accessible. Then there exists a residual subset $\mathcal{R}\subset C^0(M,\mathfrak{G})$ such that any $A\in \mathcal{R}$ 
either has all Lyapunov exponents equal at almost every point, or there exists a dominated splitting of $M \times \mathbb{R}^n$ which coincides with the Oseledets splitting almost everywhere. 
\end{lemma}

\noindent \textbf{Proof of Theorem C.} 
Since $\mathcal{S}^\ast$ is open and dense in $\mathcal{S}$ by Lemma \ref{lem:open_dense}, it suffices to find a residual subset inside $\mathcal{S}^\ast$. 
Since $Q\mathcal{S}^\ast Q$ is accessible by Lemma~\ref{main}, it follows from Lemma~\ref{corollaryBV}, that there exists a residual subset $\hat{\mathcal{R}}$ of $Q\mathcal{S}^\ast Q$ such that each normal cocycle $\hat{S} \in \hat{\mathcal{R}}$ satisfies the Bochi-Viana dichotomy. Now, with $\Theta$ defined as in (\ref{eq:Theta}), set  $\mathcal{R}:= \Theta^{-1} (\hat{\mathcal{R}} \times \hat{\Delta}_n) \cap \mathcal{S}^\ast$. Since $\Theta$ is a homeomorphism,  $\mathcal{R}$  is a  residual subset of $ \mathcal{S}^\ast$. 

It hence remains to show that the Oseledets splitting of each $S \in \mathcal{R}$ is dominated. 
By intersecting $\mathcal{R}$ with the open and dense set given by Corollary \ref{corollary:almost_rpf}, we may additionally assume that always $\dim(E^0)=1$, $E^{<0}$ is an eventually contracting bundle and that $E^0 \oplus E^{<0}$ is a dominated splitting. Note that $E^{<0}$ is never trivial since $n \geq 2$. By applying the Bochi-Viana dichotomy, the proof now reduces to the following two cases.   

\begin{enumerate}
\item If the Lyapunov spectrum of $\hat{S}$ is almost surely equal to $\{\lambda\}$ for some $\lambda <0$, then Lemma \ref{angle} in combination with $\dim(E_0)=1$ implies that the Lyapunov spectrum of $S$ is almost surely equal to $\{0,\lambda\}$. The assertion then follows from the fact that $E^0 \oplus E^{<0}$ is a dominated splitting.  
\item If the Lyapunov spectrum of $\hat{S}$ is not a singleton, then, for each non-minimal Lyapunov exponent $\lambda$,  Lemma \ref{corollaryBV} implies that $\hat{E}^\lambda \oplus \hat{E}^{<\lambda}$ is a dominated splitting. It then follows from Lemma~\ref{equivalence} that the lifted splitting also is dominated. Since also $E^0 \oplus E^{<0}$ is dominated, we then obtain that the Oseledets splitting is  a dominated splitting.
\end{enumerate}
This proves the theorem.
~ \hfill $\boxempty$

\section{Proof of Theorem B}\label{sec_proof_b}
The strategy of the proof is to conjugate $\{\mathcal{L}_x/\rho_x\}$ to a family of positive operators $\{\hat{\mathcal{L}}_x\}$ with $\hat{\mathcal{L}}_x(\mathbf{1})= \mathbf{1}$ and then employ a discretization of $\mathcal{H}_r$ in order to obtain an associated stochastic cocycle. Theorem B then follows from application of Propositions \ref{prop:continuity_1} and \ref{prop:continuity_2} and Theorem C: The continuity of $\{h_x\}$,  $\{\nu_x\}$ and $\{\rho_x\}$ allows to control the conjugations and discretizations whereas Theorem C provides the approximation by a cocycle with a dominated splitting. The relevant maps and spaces are contained in the following diagram.
\[\begin{CD}
\mathcal{H}_r    @<{\psi_x}<< 
\mathcal{H}_r @>{\pi_n^x}>> V_x @>{\hbox{\footnotesize id}}>>  V_x @>{\cong}>> \R^m\\
@VV{\mathcal{L}_x/\rho_x}V        @VV{\hat{\mathcal{L}}_x}V  @VV{\hat{\mathcal{L}}_x}V @VV{{S}_x}V @VV{S_x}V \\
\mathcal{H}_r     @<{\psi_{f(x)}}<<  \mathcal{H}_r @.  \mathcal{H}_r  @>{\pi_n^{f(x)}}>> 
 V_{f(x)} @>{\cong}>> \R^m
\end{CD}
\]
We now define the objects. With $\{h_x\}$,  $\{\nu_x\}$ and $\{\rho_x\}$ as in Propositions \ref{prop:continuity_1} and \ref{prop:continuity_2}, $d\tilde{\nu}_x := h_x d\nu_x$ and for $n \in \N$ and $x \in M$, 
\begin{eqnarray*}
\psi_x &:&  \mathcal{H}_r  \to \mathcal{H}_r, \;  g \mapsto g \cdot h_x; \quad
\hat{\mathcal{L}}_x : \mathcal{H}_r  \to \mathcal{H}_r, \;  g \mapsto  {\mathcal{L}_x(g \cdot h_x)}/{\rho_x h_{f(x)}} ;\\
\pi_n^x &:&  \mathcal{H}_r \to \mathcal{H}_r , \;  g \mapsto  {\textstyle\sum_{a \in \mathcal{W}^n} }
({\textstyle{ \tilde{\nu}_x( g \cdot \mathbf{1}_{[a]})}/{\tilde{\nu}_x(\mathbf{1}_{[a]})} }) \mathbf{1}_{[a]}; \\
V_x &:&= \pi_n^x(\mathcal{H}_r); \quad    S_x : V_x \mapsto V_{f(x)},\;  g \mapsto   \pi_n^{f(x)}(\hat{\mathcal{L}}_x(g)).
\end{eqnarray*}
Observe that 
 $\|\psi_x \|_\infty = \|h_x \|_\infty$,  $ \| \pi^x_n \|_\infty  =1$ and 
$\| \pi^x_n(g) - g \|_\infty \leq D_r(g) r^n$.
Hence, $W_x := \psi^x_n(V_x)$ is $(\|h_x\|_\infty  r^n)$-dense in $\{g \in \mathcal{H}_r: D_r(g) \leq 1\}$ in the $C^0$-topology, for all $x \in X$. Furthermore, for $g \in V_x$ and $w \in [a]$, for $a \in \mathcal{W}^n$ and $\tilde{\varphi}_x := \varphi_x + \log h_x - \log \rho_{x} - \log h_{f(x)}\circ \theta$ , we have
\[ \pi_n^{f(x)} \circ \hat{\mathcal{L}}_x(g) (w) = 
\sum_{b \in \mathcal{W}^1} \frac{\tilde{\nu}_{f(x)}( \mathbf{1}_{[a]}\cdot e^{\tilde{\varphi}_x\circ \tau_b})}{\tilde{\nu}_{f(x)}( \mathbf{1}_{[a]})} g\circ\tau_b(w).
\]
Hence, by H\"older continuity of $\varphi_x$ and $h_x$, it is well known (see, e.g., Step 3 in the proof of Theorem 4.1 in \cite{Stadlbauer:2015}) that there exists $C_{\tilde{\varphi}}>0$ such that $\exp(\tilde{\varphi}_x\circ\tau_b (v) - \tilde{\varphi}_x\circ\tau_b (w)) = 1 \pm C_{\tilde{\varphi}}r^n$ for all $v,w\in[a]$. Hence, $\tilde{\nu}_{f(x)}( \mathbf{1}_{[a]}\cdot e^{\tilde{\varphi}_x\circ \tau_b})= (1 \pm C_{\tilde{\varphi}}r^n) \tilde{\varphi}_x\circ\tau_b (w)$, which then implies that 
\[ \| \hat{\mathcal{L}}_x(g) - S_x(g)\|_\infty  \leq  C_{\tilde{\varphi}}r^n  \| \hat{\mathcal{L}}_x(|g|)\|_\infty \leq   C_{\tilde{\varphi}}r^n \| g\|_\infty, \forall g \in V_x.\]  	
By application of Theorem C, it follows that for each $\varepsilon >0$, there exists $\{A_x\}$ 
with $\sup_x \|A_x - S_x\|_\infty < \varepsilon$ and $\{A_x\}$ satisfying the above dichotomy. For 
\[B_x: W_x \to W_{f(x)}, g \mapsto \rho_x \cdot \psi_{f(x)} \circ A_x \circ \psi_{x}^{-1}(g), \] 
it follows that $\|B_x - \mathcal{L}_x\|_\infty \leq  \rho_x \|h_{f(x)}\|_\infty \|1/h_{x}\|_\infty (\varepsilon + C_{\tilde{\varphi}}r^n)$. The theorem  follows from the fact that $ \rho_x$ and $h_x$ vary continuously with respect to $x$. \hfill $\square$

\begin{remark} It is worth noting, that the proof does not provide that $W_x$ is $\varepsilon$-dense in $\{g : \|g\|_\mathcal{H} \leq 1\}$ with respect to $\|g\|_\mathcal{H}$ due to the fact that $\| \pi^x_n - {\hbox{id}} \|_\mathcal{H} \leq 2 + r^n$. Moreover, one does not obtain an approximation by an operator which is acting on $\mathcal{H}$. However, by considering the above diagram for the $n$-th iterate $\mathcal{L}_x^n$, the same estimates give rise to an approximation $B_x^{(n)}: W_x \to W_{f^n(x)} $ which can be written as a relative transfer operator with respect to a potential which is constant on cylinders of length $2n$. In particular, $B_x^{(n)}$ also acts on $\mathcal{H}$.
\end{remark}

\section*{Acknowledgements}

The first author is partially supported by National Funds through FCT - ``Funda\c{c}\~ao para a Ci\^encia e a Tecnologia'', project PEst-OE/MAT/UI0212/2011, whereas the second author has been supported in part by EU Marie-Curie IRSES Brazilian-European partnership in Dynamical Systems (FP7-PEOPLE-2012-IRSES 318999 BREUDS).


\end{document}